\documentclass[11pt]{article}

\usepackage{fullpage}
\usepackage{macros}
\usepackage{amsfonts,epsfig,graphicx}
\usepackage{amsmath,amssymb,amsthm}
\usepackage{enumerate}
\usepackage{graphics}
\usepackage[numbers]{natbib}
\usepackage{epstopdf}
\usepackage[bookmarks=true,colorlinks,citecolor=blue,urlcolor=blue]{hyperref}

\newtheorem{theorem}{Theorem}
\newtheorem{lemma}{Lemma}
\newtheorem{proposition}{Proposition}
\newtheorem{corollary}{Corollary}

\newcommand{\PREDERRSQ}[2]{\ensuremath{\frac{1}{\numobs} \|\Xmat
    (#1 - #2)\|^2_2}}

\newcommand{\xtil}{\ensuremath{\tilde{x}}}

\newcommand{\wprime}{\ensuremath{w'}}

\newcommand{\Term}{\ensuremath{T}}

\newcommand{\thetahatlam}{\ensuremath{\thetahat_\lambda}}

\long\def\comment#1{}

\usepackage{color}
\newcommand{\red}[1]{\textcolor{red}{#1}}

\newcommand{\util}{\ensuremath{\widetilde{u}}}

\newcommand{\thetait}[1]{\ensuremath{\theta^{#1}}}

\newcommand{\SsetThree}{\ensuremath{\Sset_2}}

\newcommand{\HACKETWO}{\ensuremath{\Event_0 \cap \Event_1}}

\newcommand{\REGLAM}{\ensuremath{\reg}}

\newcommand{\INITVAR}{\ensuremath{\gamma}}

\newcommand{\stepsize}{\eta}
\newcommand{\SPECBALL}{\ensuremath{\Ball_2(\stepsize; \theta^t)}}

\newcommand{\FINIT}{\ensuremath{T}}
\newcommand{\thetahack}{\ensuremath{\widehat{\theta}}}
\newcommand{\uhack}{\ensuremath{\widehat{u}}}

\newcommand{\NEWXbad}{\ensuremath{X_{\mbox{\tiny{bad}}}}}

\begin{document}

\begin{center}
  {\LARGE{\bf{Optimal prediction for sparse linear models?\\ Lower bounds for coordinate-separable M-estimators}}}

  \vspace{1cm}

  {\large
\begin{tabular}{ccccc}
Yuchen Zhang$^\star$ & & Martin J.\ Wainwright$^{\star, \dagger}$ &&
Michael I. Jordan$^{\star,\dagger}$
\end{tabular}
}

  \vspace{.5cm}

  \texttt{\{yuczhang,wainwrig,jordan\}@berkeley.edu} \\

  \vspace{.5cm}

  {\large $^\star$Department of Electrical Engineering and
  Computer Science
  ~~~~ $^\dagger$Department of Statistics} \\
\vspace{.1cm}

  {\large University of California, Berkeley} \\
  \vspace{.5cm}
\end{center}

\vspace*{.2cm}

\begin{abstract}
For the problem of high-dimensional sparse linear regression, it is
known that an $\ell_0$-based estimator can achieve a $1/n$ ``fast''
rate on the prediction error without any conditions on the design
matrix, whereas in the absence of restrictive conditions on the design
matrix, popular polynomial-time methods only guarantee the
$1/\sqrt{n}$ ``slow'' rate. In this paper, we show that the slow rate
is intrinsic to a broad class of M-estimators. In particular, for
estimators based on minimizing a least-squares cost function together
with a (possibly non-convex) coordinate-wise separable regularizer,
there is always a ``bad'' local optimum such that the associated
prediction error is lower bounded by a constant multiple of
$1/\sqrt{n}$. For convex regularizers, this lower bound applies to all
global optima.  The theory is applicable to many popular estimators,
including convex $\ell_1$-based methods as well as M-estimators based
on nonconvex regularizers, including the SCAD penalty or the MCP
regularizer.  In addition, we show that the bad local optima are very
common, in that a broad class of local minimization algorithms with
random initialization will typically converge to a bad solution.
\end{abstract}

%%%%%%%%%%%%%%%%%%%%%%%%%%%%%%%%%%%%%%%%%%%%%%%%%%%%%%%%%%%%%%%%%%%%%%%%
\section{Introduction}

The classical notion of minimax risk, which plays a central role in
decision theory, allows for the statistician to implement any possible
estimator, regardless of its computational cost.  For many problems,
there are a variety of estimators, which can be ordered in terms of
their computational complexity.  Given that it is usually feasible
only to implement polynomial-time methods, it has become increasingly
important to study computationally-constrained analogues of the
minimax estimator, in which the choice of estimator is restricted to a
subset of computationally efficient estimators~\cite{Wai14_ICM}.  A
fundamental question is when such computationally-constrained forms of
minimax risk estimation either coincide or differ in a fundamental way
from their classical counterpart.

The goal of this paper is to explore such gaps between classical and
computationally practical minimax risks, in the context of prediction
error for high-dimensional sparse regression.  Our main contribution
is to establish a fundamental gap between the classical minimax
prediction risk and the best possible risk achievable by a broad class
of $M$-estimators based on coordinate-separable regularizers, one
which includes various nonconvex regularizers that are used in
practice.

In more detail, the classical linear regression model is based on a
response vector $\yvec \in \real^\numobs$ and a design matrix $\Xmat
\in \real^{\numobs \times \usedim}$ that are linked via the relationship
\begin{align}
\label{eqn:standard-linear-model}
\yvec & = \Xmat \thetastar + w,
\end{align}
where the vector $w \in \real^\numobs$ is a random noise vector. 
Our goal is to estimate the unknown regression vector $\thetastar \in
\real^\usedim$.  Throughout this paper, we focus on the 
standard Gaussian model, in which the entries
of the noise vector $w$ are i.i.d.~$N(0, \sigma^2)$ variates, and the
case of deterministic design, in which the matrix $\Xmat$ is viewed as
non-random.  In the sparse variant of this model, the regression
vector is assumed to have a small number of non-zero coefficients. In
particular, for some positive integer $\kdim < \usedim$, the vector
$\thetastar$ is said to be $\kdim$-sparse if it has at most $\kdim$
non-zero coefficients.  Thus, the model is parameterized by the triple
$(\numobs, \usedim, \kdim)$ of sample size $\numobs$, ambient
dimension $\usedim$, and sparsity $\kdim$.  We use $\Ball_0(\kdim)$ to
the denote the $\ell_0$-``ball'' of all $\usedim$-dimensional vectors
with at most $\kdim$ non-zero entries.

An estimator $\thetahat$ is a measurable function of the pair $(y,
X)$, taking values in $\real^\usedim$, and its quality can be assessed
in different ways.  In this paper, we focus on its \emph{fixed design
  prediction error}, given by $\Exs \big[ \frac{1}{\numobs} \|X(
  \thetahat - \thetastar)\|_2^2 \big]$, a quantity that measures how
well $\thetahat$ can be used to predict the vector $X \thetastar$ of
noiseless responses. The worst-case prediction error of an estimator
$\thetahat$ over the set $\Ball_0(\kdim)$ is given by
\begin{align}
  \MSE(\thetahat; \Xmat) & \defn \sup_{\thetastar\in \Ball_0(\kdim)}
  \frac{1}{\numobs}\E[\ltwos{\Xmat (\thetahat - \thetastar)}^2]
\end{align}
Given that $\thetastar$ is $\kdim$-sparse, the most direct approach
would be to seek a $\kdim$-sparse minimizer to the least-squares cost
$\|\yvec - \Xmat \theta\|_2^2$, thereby obtaining the $\ell_0$-based
estimator
\begin{align}
\label{EqnDefnEllZeroEstimator}
\thetazero & \in \arg \min_{\theta\in \Ball_0(\kdim)} \|\yvec - \Xmat
\theta\|_2^2.
\end{align}
The $\ell_0$-based estimator $\thetazero$ is
known~\cite{BunWegTsyb07,raskutti2011minimax} to satisfy a bound of
the form
\begin{align}
\label{eqn:l0-optimal-rate}
  \MSE(\thetazero; \Xmat) \precsim \frac{\sigma^2 \, \kdim \log
    \usedim}{\numobs},
\end{align}
where $\precsim$ denotes an inequality up to constant factors
(independent of the triple $(\numobs, \usedim, \kdim)$ as well as the
standard deviation $\sigma$).  However, it is not tractable to compute
this estimator in a brute force manner, since there are ${\usedim
  \choose \kdim}$ subsets of size $\kdim$ to consider.

The computational intractability of the $\ell_0$-based estimator has
motivated the use of various heuristic algorithms and approximations,
including the basis pursuit method~\cite{chen1998atomic}, the Dantzig
selector~\cite{candes2007dantzig}, as well as the extended family of
Lasso estimators~\cite{tibshirani1996regression,
  chen1998atomic,zou2006adaptive,belloni2011square}.  Essentially,
these methods are based on replacing the $\ell_0$-constraint with its
$\ell_1$-equivalent, in either a constrained or penalized form.  There
is now a very large body of work on the performance of such methods,
covering different criteria including support recovery, $\ell_2$-norm
error and prediction error (e.g., see the book~\cite{BuhVan11} and
references therein).

For the case of fixed design prediction error that is the primary
focus here, such $\ell_1$-based estimators are known to achieve the
bound~\eqref{eqn:l0-optimal-rate} only if the design matrix $\Xmat$
satisfies certain conditions, such as the restricted eigenvalue (RE)
condition or compatibility
condition~\cite{bickel2009simultaneous,van2009conditions} or the
stronger restricted isometry property~\cite{candes2007dantzig}; see
the paper~\cite{van2009conditions} for an overview of these various
conditions, and their inter-relationships. Without such conditions,
the best known guarantees for $\ell_1$-based estimators are of the
form
\begin{align}
\label{eqn:l1-achievable-rate}
  \MSE(\thetahat_{\ell_1}; \Xmat) \precsim \sigma \, R \,
  \sqrt{\frac{\log \usedim}{\numobs}},
\end{align}
a bound that is valid without any RE conditions on the design matrix
$\Xmat$ whenever the $\kdim$-sparse regression vector $\thetastar$ has
$\ell_1$-norm bounded by $R$ (e.g., see the
papers~\cite{BunWegTsyb07,Nem00,raskutti2011minimax}.)

The substantial gap between the ``fast''
rate~\eqref{eqn:l0-optimal-rate} and the ``slow''
rate~\eqref{eqn:l1-achievable-rate} leaves open a fundamental
question: is there a computationally efficient estimator attaining the
bound~\eqref{eqn:l0-optimal-rate} for general design matrices?  In the
following subsections, we provide an overview of the currently known
results on this gap, and we then provide a high-level statement of
the main result of this paper.

%%%%%%%%%%%%%%%%%%%%%%%%%%%%%%%%%%%%%%%%%%%%%%%%%%%%%%%%%%%%%%%%%%%%%%%

\subsection{Lower bounds for Lasso}

Given the gap between the fast rate~\eqref{eqn:l0-optimal-rate} and
the Lasso's slower rate~\eqref{eqn:l1-achievable-rate}, one
possibility might be that existing analyses of prediction error are
overly conservative, and $\ell_1$-based methods can actually achieve
the bound~\eqref{eqn:l0-optimal-rate}, without additional constraints
on $\Xmat$.  Some past work has given negative answers to this
quesiton.  Foygel and Srebro~\cite{FoySre11} constructed a 2-sparse
regression vector and a random design matrix for which the Lasso
prediction error with any choice of regularization parameter
$\regparn$ is lower bounded by $1/\sqrt{\numobs}$. In particular,
their proposed regression vector is $\thetastar
=(0,\dots,0,\frac{1}{2},\frac{1}{2})$.  In their design matrix, the
columns are randomly generated with distinct covariances, and
moreover, such that the rightmost column is strongly correlated with
the other two columns on its left. With this particular regression
vector and design matrix, they show that Lasso's prediction error is
lower bounded by $1/\sqrt{\numobs}$ for \emph{any} choice of Lasso
regularization parameter $\lambda$. This construction is explicit for
Lasso, and thus does not apply to more general M-estimators. Moreover,
for this particular counterexample, there is a one-to-one
correspondence between the regression vector and the design matrix, so
that one can identify the non-zero coordinates of $\thetastar$ by
examining the design matrix. Consequently, for this construction, a
simple reweighted form of the Lasso can be used to achieve the fast
rate.  In particular, the reweighted Lasso estimator
\begin{align}
\label{eqn:general-reweighed-lasso}
\thetahat_{\rm wl} \in \arg \min_{\theta\in \R^d} \left \{ \|\yvec -
\Xmat \theta\|_2^2 + \lambda \sum_{j=1}^d \alpha_j |\theta_j| \right
\},
\end{align}
with $\lambda$ chosen in the usual manner ($\lambda \asymp \sigma
\sqrt{\frac{\log \usedim}{\numobs}}$), weights $\alpha_{d-1} =
\alpha_d = 1$, and the remaining weights $\{\alpha_1, \ldots,
\alpha_{d-2} \}$ chosen to be sufficiently large, has this property.
Dalalyan et al.~\cite{dalalyan2014prediction} construct a stronger
counter-example, for which the prediction error of Lasso is again
lower bounded by $1/\sqrt{n}$.  For this counterexample, there is no
obvious correspondence between the regression vector and the design
matrix. Nevertheless, as we show in
Appendix~\ref{sec:fast-rate-dalalyan}, the reweighted Lasso
estimator~\eqref{eqn:general-reweighed-lasso} with a proper choice of
the regularization coefficients still achieves the fast rate on this
example.  Another related piece of work is by Cand{\`e}s and
Plan~\cite{candes2009near}. They construct a design matrix for which
the Lasso estimator, when applied with the usual choice of
regularization parameter $\lambda \asymp \sigma (\frac{\log
    \usedim}{\numobs})^{1/2}$, has sub-optimal prediction error.  Their
matrix construction is spiritually similar to ours, but the
theoretical analysis is limited to the Lasso for a particular choice
of regularization parameter. It does not rule out the possibility that
other choices of regularization parameters, or other polynomial-time
estimators can achieve the fast rate. In contrast, our hardness result
applies to general $M$-estimators based on coordinatewise separable
regularizers, and it allows for arbitrary regularization parameters.

%%%%%%%%%%%%%%%%%%%%%%%%%%%%%%%%%%%%%%%%%%%%%%%%%%%%%%%%%%%%%%%%%%%%%%%%

\subsection{Complexity-theoretic lower bound for polynomial-time sparse estimators}

In our own recent work~\cite{zhang2014lower}, we have provided a
complexity-theoretic lower bound that applies to a very broad class of
polynomial-time estimators.  The analysis is performed under a
standard complexity-theoretic condition---namely, that the class $\np$
is not a subset of the class $\ppoly$---and shows that there is no
polynomial-time algorithm that returns a $\kdim$-sparse vector that
achieves the fast rate. The lower bound is established as a function
of the restricted eigenvalue of the design matrix. Given sufficiently
large $(\numobs,\kdim,\usedim)$ and any $\gamma > 0$, a design matrix
$X$ with restricted eigenvalue $\gamma$ can be constructed, such that
every polynomial-time $k$-sparse estimator $\thetahat_{\rm poly}$ has
its minimax prediction risk lower bounded as
\begin{align}
\label{eqn:polytime-paper-lower-bound}
 \MSE(\thetahat_{\rm poly}; \Xmat) \succsim \frac{\sigma^2 \,
   \kdim^{1-\delta} \log \usedim}{\gamma \numobs},
\end{align}
where $\delta > 0$ is an arbitrarily small positive scalar. Note that
the fraction $\kdim^{-\delta}/\gamma$, which characterizes the gap between the fast rate and the rate~\eqref{eqn:polytime-paper-lower-bound}, could be arbitrarily large. The lower bound has the following consequence:
any estimator that achieves the fast rate must either not be polynomial-time,
or must return a regression vector that is not $k$-sparse.

The condition that the estimator is $k$-sparse is essential in the
proof of lower bound~\eqref{eqn:polytime-paper-lower-bound}. In
particular, the proof relies on a reduction between estimators with
low prediction error in the sparse linear regression model, and
methods that can solve the 3-set covering
problem~\cite{natarajan1995sparse}, a classical problem that is known
to be NP-hard. The 3-set covering problem takes as input a list of
3-sets, which are subsets of a set $\mathcal{S}$ whose cardinality is
$3k$. The goal is to choose $k$ of these subsets in order to cover the
set $\mathcal{S}$. The lower
bound~\eqref{eqn:polytime-paper-lower-bound} is established by showing
that if there is a $k$-sparse estimator achieving better prediction
error, then it provides a solution to the 3-set covering problem, as
every non-zero coordinate of the estimate corresponds to a chosen
subset. This hardness result does not eliminate the possibility of
finding a polynomial-time estimator that returns dense vectors
satisfying the fast rate. In particular, it is possible that a dense
estimator cannot be used to recover a a good solution to the 3-set
covering problem, implying that it is not possible to use the hardness
of $3$-set covering to assert the hardness of achieving low prediction
error in sparse regression.

At the same time, there is some evidence that better prediction error
can be achieved by dense estimators. For instance, suppose that we
consider a sequence of high-dimensional sparse linear regression
problems, such that the restricted eigenvalue $\gamma =
\gamma_\numobs$ of the design matrix $X \in \real^{\numobs \times
  \usedim}$ decays to zero at the rate $\gamma_\numobs = 1/\numobs^2$.
For such a sequence of problems, as $\numobs$ diverges to infinity,
the lower bound~\eqref{eqn:polytime-paper-lower-bound}, which applies
to $\kdim$-sparse estimators, goes to infinity, whereas the Lasso
upper bound~\eqref{eqn:l1-achievable-rate} converges to zero.
Although this behavior is somewhat mysterious, it is not a
contradiction.  Indeed, what makes Lasso's performance better than the
lower bound~\eqref{eqn:polytime-paper-lower-bound} is that it allows
for non-sparse estimates. In this example, truncating the Lasso's
estimate to be $k$-sparse will substantially hurt the prediction
error.  In this way, we see that proving lower bounds for non-sparse
estimators---the problem to be addressed in this paper---is a
substantially more challenging task than proving lower bound for estimators
that must return sparse outputs.

%%%%%%%%%%%%%%%%%%%%%%%%%%%%%%%%%%%%%%%%%%%%%%%%%%%%%%%%%%%%%%%%%%%%%%%%

\subsection{Main results of this paper}

With this context in place, let us now turn to a high-level statement
of the main results of this paper.  More precisely, our contribution
is to provide additional evidence against the polynomial achievability
of the fast rate~\eqref{eqn:l0-optimal-rate}, in particular by showing
that the slow rate~\eqref{eqn:l1-achievable-rate} is a lower bound for
a broad class of M-estimators, namely those based on minimizing a
least-squares cost function together with a coordinate-wise
decomposable regularizer.  In particular, we consider estimators that
are based on an objective function of the form $L(\theta; \lambda) =
\frac{1}{\numobs} \|y - X \theta\|_2^2 + \lambda \, \REGLAM(\theta)$,
for a weighted regularizer $\REGLAM: \real^\usedim \rightarrow \real$
that is coordinate-separable.  See Section~\ref{SecCoord} for a
precise definition of this class of estimators.  Our first main result
(Theorem~\ref{theorem:main-lower-bound}) establishes that there is
always a matrix $\Xmat\in \R^{\numobs\times \usedim}$ such that for
any coordinate-wise separable function~$\reg$ and for any choice of
weight $\lambda \geq 0$, the objective $L$ always has at least one
local optimum $\thetahat_\lambda$ such that
\begin{align}
\label{eqn:intro-main-result}
  \sup_{\thetastar\in \Ball_0(\kdim)} \E \Big[ \frac{1}{\numobs}
    \ltwos{\Xmat (\thetahat_\lambda - \thetastar \big)}^2 \Big]
  \succsim \frac{\sigma}{\sqrt{\numobs}}.
\end{align}
Moreover, if the regularizer $\REGLAM$ is convex, then this lower
bound applies to all global optima of the convex criterion $L$.  This
lower bound is applicable to many popular estimators, including the
ridge regression estimator~\cite{hoerl1970ridge}, the basis pursuit
method~\cite{chen1998atomic}, the Lasso
estimator~\cite{tibshirani1996regression}, the weighted Lasso
estimator~\cite{zou2006adaptive}, the square-root Lasso
estimator~\cite{belloni2011square}, and least squares based on
nonconvex regularizers such as the SCAD penalty~\cite{fan2001variable}
or the MCP penalty~\cite{zhang2010nearly}.

In the nonconvex setting, it is impossible (in general) to guarantee
anything beyond local optimality for any solution found by a
polynomial-time algorithm~\cite{ge2015strong}. Nevertheless, to play
the devil's advocate, one might argue that the assumption that an
adversary is allowed to pick a bad local optimum could be overly
pessimistic for statistical problems.  In order to address this
concern, we prove a second result
(Theorem~\ref{theorem:gradient-descent-lower-bound}) that demonstrates
that bad local solutions are difficult to avoid.  Focusing on a class
of local descent methods, we show that given a random isotropic
initialization centered at the origin, the resulting stationary points
have poor mean-squared error---that is, they can only achieve the slow
rate.  In this way, this paper shows that the gap between the fast and
slow rates in high-dimensional sparse regression cannot be closed via
standard application of a very broad class of methods.  In conjunction
with our earlier complexity-theoretic paper~\cite{zhang2014lower}, it
adds further weight to the conjecture that there is a fundamental gap
between the performance of polynomial-time and exponential-time
methods for sparse prediction.

The remainder of this paper is organized as follows.  We begin in
Section~\ref{SecBackground} with further background, including a
precise definition of the family of $M$-estimators considered in this
paper, some illustrative examples, and discussion of the prediction
error bound achieved by the Lasso.  Section~\ref{sec:main-result} is
devoted to the statements of our main results, along with discussion
of their consequences.  In Section~\ref{SecProofs}, we provide the
proofs of our main results, with some technical lemmas deferred to the
appendices.  We conclude with a discussion in
Section~\ref{SecDiscussion}.

%%%%%%%%%%%%%%%%%%%%%%%%%%%%%%%%%%%%%%%%%%%%%%%%%%%%%%%%%%%%%%%%%%%%%%%%%%%

\section{Background and problem set-up}
\label{SecBackground}

As previously described, an instance of the sparse linear regression
problem is based on observing a pair $(\Xmat, \yvec)\in
\R^{\numobs\times \usedim} \times \R^\numobs$ of instances that are
linked via the linear model~\eqref{eqn:standard-linear-model}, where
the unknown regressor $\thetastar$ is assumed to be $\kdim$-sparse,
and so belongs to the $\ell_0$-ball $\Ball_0(\kdim)$.  Our goal is to
find a good predictor, meaning a vector $\thetahat$ such that the
mean-squared prediction error $\frac{1}{\numobs} \ltwos{\Xmat(
  \thetahat -\thetastar)}^2$ is small.

\subsection{Least squares with coordinate-separable regularizers}
\label{SecCoord}

The analysis of this paper applies to estimators that are based on
minimizing a cost function of the form
\begin{align}
\label{EqnLoss}
L(\theta; \lambda) & = \frac{1}{\numobs} \|y - X \theta\|_2^2 +
\lambda \, \REGLAM(\theta),
\end{align}
where $\REGLAM:\real^\usedim \rightarrow \real$ is a
\emph{regularizer}, and $\lambda \geq 0$ is a regularization weight.
We consider the following family $\family$ of coordinate-separable
regularizers:
\begin{enumerate}[(i)]
\item The function $\REGLAM: \real^\usedim \rightarrow \real$ is
  coordinate-wise decomposable, meaning that $\REGLAM(\theta) =
  \sum_{j=1}^\usedim \reg_j(\theta_j)$ for some univariate functions
  $\rho_j:\real \rightarrow \real$.
 \item Each univariate function satisfies $\rho_j(0)=0$ and is
   symmetric around zero (i.e., $\rho_j(t) = \rho_j(-t)$ for all $t\in
   \R$).
  \item On the nonnegative real line $[0,+\infty)$, each function
    $\rho_j$ is nondecreasing.
\end{enumerate}

\noindent Let us consider some examples to illustrate this definition.\\

\paragraph{Bridge regression:} The family of bridge regression estimates~\cite{frank1993statistical}
take the form
\begin{align*}
  \thetahat_{\tiny{\mbox{bidge}}} \in \arg \min_{ \theta \in
    \real^\usedim} \Big \{ \frac{1}{\numobs} \|\yvec - \Xmat
  \theta\|_2^2 + \lambda \sum_{i=1}^\usedim |\theta|^\gamma \Big \}.
\end{align*}
Note that this is a special case of the objective
function~\eqref{EqnLoss} with $\rho_j(\cdot) = |\cdot|^\gamma$ for
each coordinate.  When $\gamma \in \{1,2\}$, it corresponds to the
Lasso estimator and the ridge regression estimator respectively.  The
analysis of this paper provides lower bounds for both estimators,
uniformly over the choice of $\lambda$.

\paragraph{Weighted Lasso:} The weighted Lasso estimator~\cite{zou2006adaptive}
 uses a weighted $\ell_1$-norm to regularize the empirical risk, and
 leads to the estimator
\begin{align*}
\thetahat_{\tiny{\mbox{wl}}} \in \arg \min_{ \theta \in \real^\usedim}
\Big \{ \frac{1}{\numobs} \|\yvec - \Xmat \theta\|_2^2 + \lambda
\sum_{i=1}^\usedim \alpha_i|\theta_i| \Big \}.
\end{align*}
Here $\alpha_1, \dots, \alpha_\usedim$ are weights that can be
adaptively chosen with respect to the design matrix $\Xmat$.  The
weighted Lasso can perform better than the ordinary Lasso,
corresponding to the special case in which all $\alpha_j$ are all
equal.  For instance, on the counter-example proposed by Foygel and
Srebro~\cite{FoySre11}, for which the ordinary Lasso estimator
achieves only the slow $1/\sqrt{\numobs}$ rate, the weighted Lasso
estimator achieves the $1/\numobs$ convergence rate.  Nonetheless, the
analysis of this paper shows that there are design matrices for which
the weighted Lasso, even when the weights are chosen adaptively with
respect to the design, has prediction error at least a constant
multiple of $1/\sqrt{\numobs}$.

\paragraph{Square-root Lasso:} The square-root Lasso 
estimator~\cite{belloni2011square} is defined by minimizing the
criterion
\begin{align*}
  \thetahat_{\tiny{\mbox{sqrt}}} \in \arg \min_{ \theta \in
    \real^\usedim} \Big \{ \frac{1}{\sqrt{\numobs}} \|\yvec - \Xmat
  \theta\|_2 + \lambda \lone{\theta} \Big \}.
\end{align*}
This criterion is slightly different from our general objective
function~\eqref{EqnLoss}, since it involves the square root of the
least-squares error. Relative to the Lasso, its primary advantage is
that the optimal setting of the regularization parameter does not
require the knowledge of the standard deviation of the noise.  For the
purposes of the current analysis, it suffices to note that by
Lagrangian duality, every square-root Lasso estimate
$\thetahat_{\tiny{\mbox{sqrt}}}$ is a minimizer of the least-squares
criterion $\|\yvec - \Xmat \theta\|_2$, subject to $\|\theta\|_1 \leq R$,
for some radius $R \geq 0$ depending on $\lambda$.  Consequently, as
the weight $\lambda$ is varied over the interval $[0,\infty)$, the
  square root Lasso yields the same solution path as the Lasso.  Since
  our lower bounds apply to the Lasso for any choice of $\lambda \geq
  0$, they also apply to all square-root Lasso solutions.

\begin{figure}
\centering
\includegraphics[width = 0.65\textwidth]{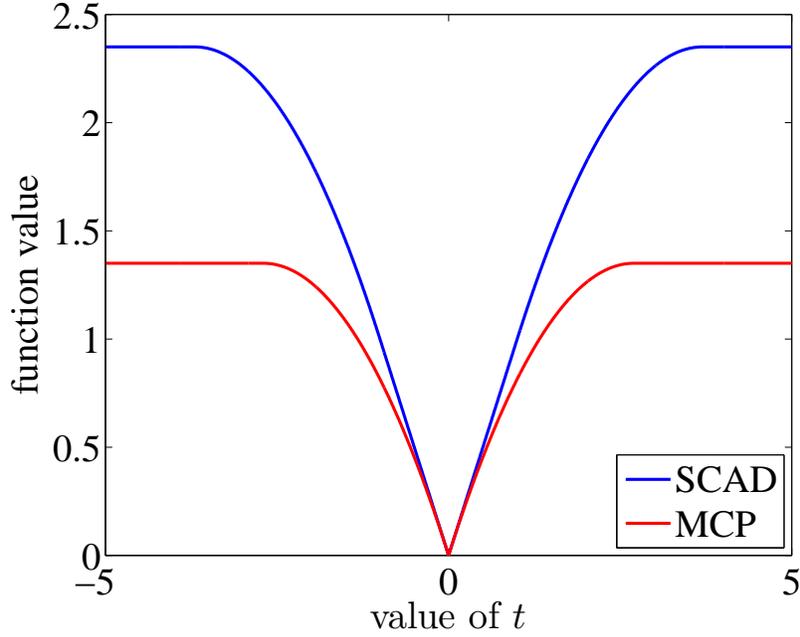}
\caption{Plots with regularization weight $\lambda = 1$, and
  parameters $a = 3.7$ for SCAD, and $b = 2.7$ for MCP.}
\label{fig:penalty-func}
\end{figure}

\paragraph{SCAD penalty or MCP regularizer:} 
Due to the intrinsic bias induced by $\ell_1$-regularization, various
forms of nonconvex regularization are widely used.  Two of the most
popular are the SCAD penalty, due to Fan and
Li~\cite{fan2001variable}, and the MCP penalty, due to
Zhang et al.~\cite{zhang2010nearly}.  The family of SCAD penalties takes the
form
\begin{align*}
\phi_\lambda(t) \defeq \frac{1}{\lambda} \begin{cases} \lambda |t| &
  \mbox{for $|t|\leq \lambda$},\\ -(t^2-2 a \lambda |t| +
  \lambda^2)/(2a-2) & \mbox{for $\lambda < |t| \leq
    a\lambda$},\\ (a+1)\lambda^2/2 & \mbox{for $|t| \geq a\lambda$},
  \end{cases}
\end{align*}
where $a > 2$ is a fixed parameter.  When used with the least-squares
objective, it is a special case of our general set-up with
$\reg_j(\theta_j) = \phi_\lambda(\theta_j)$ for each coordinate $j =
1, \ldots, \usedim$.  Similarly, the MCP penalty takes the form
\begin{align*}
  \reg_\lambda(t) \defeq \int_0^{|t|} \left( 1 - \frac{z}{\lambda
    b}\right)_+ d z,
\end{align*}
where $b > 0$ is a fixed parameter. It can be verified that both the
SCAD penalty and the MCP regularizer belong to the function class
$\family$ previously defined. See Figure~\ref{fig:penalty-func} for
a graphical illustration of the SCAD penalty and the MCP regularizer.

%%%%%%%%%%%%%%%%%%%%%%%%%%%%%%%%%%%%%%%%%%%%%%%%%%%%%%%%%%%%%%%%%%%%%%%

\subsection{Prediction error for the Lasso}

We now turn to a precise statement of the best known upper bounds for
the Lasso prediction error. We assume that the design matrix satisfies
the column normalization condition. More precisely, letting $X_j \in
\real^\numobs$ denote the $j^{th}$ column of the design matrix
$\Xmat$, we say that it is $1$-column normalized if
\begin{align}
\label{eqn:column-normalization-condition}
\frac{\ltwos{X_j}}{\sqrt{n}} \leq 1 \qquad \mbox{for $j = 1,2,
  \ldots,\usedim$.}
\end{align}
Our choice of the constant $1$ is to simplify notation; the more
general notion allows for an arbitrary constant $C$ in this bound.

In addition to the column normalization condition, if the design matrix further satisfies a restricted eigenvalue (RE)
condition~\cite{bickel2009simultaneous,van2009conditions}, then the
Lasso is known to achieve the fast rate~\eqref{eqn:l0-optimal-rate}
for prediction error. More precisely, restricted eigenvalues are defined
in terms of subsets $\PlainSset$ of the index set $\{1, 2, \ldots,
\usedim \}$, and a cone associated with any such subset.  In
particular, letting $\PlainSbar$ denote the complement of
$\PlainSset$, we define the cone
\begin{align*}
%\label{EqnDefnConeSet}
\ConeSet(\PlainSset) & \defn \big \{ \theta \in \real^\usedim \, \mid
\, \|\theta_{\PlainSbar}\|_1 \leq 3 \|\theta_{\PlainSset}\|_1 \big \}.
\end{align*}
Here $\|\theta_{\PlainSbar}\|_1 \defn \sum_{j \in \PlainSbar}
|\theta_j|$ corresponds to the $\ell_1$-norm of the coefficients
indexed by $\PlainSbar$, with $\|\theta_{\PlainSset}\|_1$ defined
similarly.  Note that any vector $\thetastar$ supported on
$\PlainSset$ belongs to the cone $\ConeSet(\PlainSset)$; in addition,
it includes vectors whose $\ell_1$-norm on the ``bad'' set
$\PlainSbar$ is small relative to their $\ell_1$-norm on $\PlainSset$.
Given triplet $(\numobs, \usedim, \kdim)$, the matrix $\Xmat \in
\real^{\numobs \times \usedim}$ is said to satisfy a $\RECON$-RE
condition (also known as a compatibility condition) if
\begin{align}
\label{EqnDefnRE}
\frac{1}{\numobs} \|X \theta\|_2^2 & \geq \RECON \|\theta\|_2^2 \qquad
\mbox{for all $\theta \in \bigcup
\limits_{|\PlainSset| = \kdim} \ConeSet(\PlainSset)$.}
\end{align}
The following
result~\cite{bickel2009simultaneous,negahban2012unified,BuhVan11}
provides a bound on the prediction error for the Lasso estimator:

\begin{proposition}[Prediction error for Lasso with RE condition]
\label{PropLassoThresh}
Consider the standard linear model for a design matrix $\Xmat$
satisfying the column normalization
condition~\eqref{eqn:column-normalization-condition} and the
$\gamma$-RE condition. Then for any vector $\thetastar\in
\mathbb{B}_0(k)$, the Lasso estimator $\thetahat_{\lambda_n}$ with
$\regparn = 4 \sigma \sqrt{\frac{\log \usedim}{\numobs}}$ satisfies
\begin{align}
\label{EqnLassoThreshBound}
\frac{1}{\numobs} \| \Xmat \thetahat_{\lambda_n} - \Xmat
\thetastar\|_2^2 & \leq \frac{64\sigma^2 \kdim \log
  \usedim}{\gamma^2\numobs} \qquad \mbox{for any $\thetastar \in
  \Ball_0(\kdim)$}
\end{align}
with probability at least $1 - c_1 e^{-c_2 n \lambda_n^2}$.
\end{proposition}

The Lasso rate~\eqref{EqnLassoThreshBound} will match the optimal
rate~\eqref{eqn:l0-optimal-rate} if the RE constant~$\gamma$ is
bounded away from zero. If $\gamma$ is close to zero, then the Lasso
rate could be arbitrarily worse than the optimal rate.  It is known
that the RE condition is necessary for recovering the true vector
$\thetastar$~\cite[see, e.g.][]{raskutti2011minimax}, but minimizing
the prediction error should be easier than recovering the true vector.
In particular, strong correlations between the columns of $X$, which
lead to violations of the RE conditions, should have no effect on the
intrinsic difficulty of the prediction problem.  Recall that the
$\ell_0$-based estimator $\thetazero$ satisfies the prediction error
upper bound~\eqref{eqn:l0-optimal-rate} without any constraint on the
design matrix.  Moreover, Raskutti et al.~\cite{raskutti2011minimax}
show that many problems with strongly correlated columns are actually
easy from the prediction point of view.

In the absence of RE conditions, $\ell_1$-based methods are known to
achieve the slow $1/\sqrt{\numobs}$ rate, with the only constraint on
the design matrix being a uniform column
bound~\cite{bickel2009simultaneous}:

\begin{proposition}[Prediction error for Lasso without RE condition]
\label{PropLasso}
Consider the standard linear model for a design matrix $\Xmat$
satisfying the column normalization
condition~\eqref{eqn:column-normalization-condition}.  Then for any
vector $\thetastar\in \Ball_0(\kdim)\cap \Ball_1(\radius)$, the Lasso
estimator $\thetaregparn$, with $\regparn = 4 \sigma \sqrt{\frac{2\log
    \usedim}{\numobs}}$, satisfies the bound
\begin{align}
\label{EqnLassoBound}
\PREDERRSQ{\thetaregparn}{\thetastar} & \leq \UNICON\; \sigma \radius
\Big ( \sqrt{\frac{2\log \usedim}{\numobs}} + \delta \Big),
\end{align}
with probability at least $1 - c_1 \usedim \, e^{-c_2 \numobs
  \delta^2}$.
\end{proposition}
\noindent Combining the bounds of
Proposition~\ref{PropLassoThresh} and Proposition~\ref{PropLasso}, we
have
\begin{align}\label{eqn:combined-upper-bound}
  \MSE(\thetahat_{\ell_1}; \Xmat) \leq \UNICON'
  \min\Big\{\frac{\sigma^2 k \log d}{\gamma^2 n}, \sigma \radius
  \sqrt{\frac{ \log \usedim}{\numobs}}\Big\}.
\end{align}
If the RE constant $\gamma$ is sufficiently close to zero, then the
second term on the right-hand side will dominate the first term. In
that case, the $1/\sqrt{\numobs}$ achievable rate is substantially
slower than the $1/\numobs$ optimal rate for reasonable ranges of
$(\kdim, R)$. One might wonder whether the analysis leading to the
bound~\eqref{eqn:combined-upper-bound} could be sharpened so as to
obtain the fast rate.  Among other consequences, our first main result
(Theorem~\ref{theorem:main-lower-bound} below) shows that no
substantial sharpening is possible.

%%%%%%%%%%%%%%%%%%%%%%%%%%%%%%%%%%%%%%%%%%%%%%%%%%%%%%%%%%%%%%%%%%%%%%%%%%

\section{Main results}
\label{sec:main-result}

We now turn to statements of our main results, and discussion of their
consequences.

\subsection{A general lower bound}

Our analysis applies to the set of local minima of the objective
function $L$ defined in equation~\eqref{EqnLoss}.  More precisely, a
vector $\thetatil$ is a local minimum of the function $\theta \mapsto
L(\theta;\lambda)$ if there is an open ball $\Ball$ centered at
$\thetatil$ such that $\thetatil \in \arg \min \limits_{\theta\in
  \Ball} L(\theta;\lambda)$.  We then define the set
\begin{align}
\label{eqn:define-solution-path}
  \Thetahat_\lambda \defeq \Big \{ \theta \in \real^\usedim \, \mid \,
  \mbox{$\theta$ is a local minimum of the function $\theta \mapsto
    L(\theta;\lambda)$} \Big\},
\end{align}
an object that depends on the triplet $(\Xmat,\yvec,\rho)$ as well as
the choice of regularization weight $\lambda$.  Since the function
$\reg$ might be non-convex, the set $\Thetahat_\lambda$ may contain
multiple elements.  

At best, a typical descent method applied to the objective $L$ can be
guaranteed to converge to some element of $\Thetahat_\lambda$.  The
following theorem provides a lower bound, applicable to any method
that always returns some local minimum of the objective
function~\eqref{EqnLoss}.

\begin{theorem}
\label{theorem:main-lower-bound}
For any pair $(\numobs, \usedim)$ such that \mbox{$\usedim \geq
  \numobs \geq 4$,} any sparsity level \mbox{$\kdim \geq 2$} and any
radius \mbox{$\radius \geq 8{\sigma}/{\sqrt{\numobs}}$,} there is a
design matrix $\Xmat \in \R^{\numobs\times \usedim}$ satisfying the
column normalization
condition~\eqref{eqn:column-normalization-condition} such that for any
coordinate-separable penalty, we have
\begin{subequations}
\begin{align}
\label{eqn:main-lower-bound}
\sup_{\thetastar\in \Ball_0(\kdim)\cap \Ball_1(\radius)} \E \left[
  \inf_{\lambda\geq 0} \sup_{\theta\in \Thetahat_\lambda}
  \PREDERRSQ{\theta}{\thetastar} \right] & \geq \UNICON \min\left\{
\sigma^2, \frac{\sigma\radius}{\sqrt{\numobs}}\right\}.
\end{align}
Moreover, for any convex coordinate-separable penalty, we have
\begin{align}
\label{eqn:main-lower-bound-for-convex-penalty}
 \sup_{\thetastar\in \Ball_0(\kdim)\cap \Ball_1(\radius)} \E \left[
   \inf_{\lambda\geq 0} \inf_{\theta\in \Thetahat_\lambda}
   \PREDERRSQ{\theta}{\thetastar} \right] \geq \UNICON \min\left\{
 \sigma^2, \frac{\sigma\radius}{\sqrt{\numobs}}\right\}.
\end{align}
\end{subequations}
\end{theorem}
\noindent In both of these statements, the constant $c$ is universal,
independent of $(\numobs, \usedim, \kdim, \sigma, \radius)$ as well as
the design matrix.  See Section~\ref{sec:prove-main-result} for the
proof.

In order to interpret the lower bound~\eqref{eqn:main-lower-bound},
consider any estimator $\thetahat$ that takes values in the set
$\Thetahat_\lambda$, corresponding to local minima of $L$.  The result
is of a game-theoretic flavor: the statistician is allowed to
adaptively choose $\lambda$ based on the observations $(y, X)$,
whereas nature is allowed to act adversarially in choosing a local
minimum for every execution of $\thetahat_\lambda$.  Under this
setting, Theorem~\ref{theorem:main-lower-bound} implies that
\begin{align}
 \label{eqn:concretized-lower-bound}
\sup_{\thetastar\in \Ball_0(\kdim)\cap \Ball_1(\radius)}
\frac{1}{\numobs} \E \left[ \ltwos{\Xmat\thetahat_\lambda -
    \Xmat\thetastar}^2 \right] &\geq \UNICON \min\left\{ \sigma^2,
\frac{\sigma\radius}{\sqrt{\numobs}}\right\}.
\end{align}
For any convex regularizer (such as the $\ell_1$-penalty underlying
the Lasso estimate),
equation~\eqref{eqn:main-lower-bound-for-convex-penalty} provides a
stronger lower bound, one that holds uniformly over all choices of
$\lambda \geq 0$ and all (global) minima.  For the Lasso estimator,
the lower bound of Theorem~\ref{theorem:main-lower-bound} matches the
upper bound~\eqref{EqnLassoBound} up to the logarithmic term
$\sqrt{\log \usedim}$, showing that the lower bound is almost tight.\\

It is possible that lower bounds of this form hold only for extremely
ill-conditioned design matrices, which would render the consequences
of the result less broadly applicable.  In particular, it is natural
to wonder whether it is also possible to prove a non-trivial lower
bound when the restricted eigenvalues are bounded above zero. Recall
that under the RE condition with a positive constant~$\gamma$, the
Lasso will achieve a mixture rate~\eqref{eqn:combined-upper-bound},
consisting of a scaled fast rate $1/(\gamma^2 n)$ and the slow rate
$1/\sqrt{n}$. The following result shows that this mixture rate cannot
be improved to match the fast rate.

\begin{corollary}
\label{coro:lower-bound}
For any sparsity level $\kdim \geq 1$, any integers $\usedim = \numobs
\geq 4k$, any radius \mbox{$\radius \geq 8{\sigma}/{\sqrt{\numobs}}$}
and any constant $\gamma \in (0,1]$, there is a design matrix $\Xmat
\in \R^{\numobs\times \usedim}$ satisfying the column normalization
condition~\eqref{eqn:column-normalization-condition} and the
$\gamma$-RE condition, such that for any coordinate-separable penalty,
we have
\begin{subequations}
\begin{align}
\label{eqn:coro-lower-bound}
\sup_{\thetastar\in \Ball_0(2\kdim)\cap \Ball_1(\radius)} \E \left[
  \inf_{\lambda\geq 0} \sup_{\theta\in \Thetahat_\lambda}
  \PREDERRSQ{\theta}{\thetastar} \right] & \geq \UNICON \min\left\{
\sigma^2, \frac{k \sigma^2}{\gamma n},
\frac{\sigma\radius}{\sqrt{\numobs}}\right\}.
\end{align}
Moreover, for any convex coordinate-separable penalty, we have
\begin{align}
\label{eqn:coro-lower-bound-for-convex-penalty}
 \sup_{\thetastar\in \Ball_0(2\kdim)\cap \Ball_1(\radius)} \E \left[
   \inf_{\lambda\geq 0} \inf_{\theta\in \Thetahat_\lambda}
   \PREDERRSQ{\theta}{\thetastar} \right] \geq \UNICON \min\left\{
 \sigma^2, \frac{k \sigma^2}{\gamma n},
 \frac{\sigma\radius}{\sqrt{\numobs}}\right\}.
\end{align}
\end{subequations}
\end{corollary}

Since none of the three terms on the right-hand side of
inequalities~\eqref{eqn:coro-lower-bound}
and~\eqref{eqn:coro-lower-bound-for-convex-penalty} matches the
optimal rate~\eqref{eqn:l0-optimal-rate}, the corollary implies that
the optimal rate is not achievable even if the restricted eigenvalues
are bounded above zero. Comparing this lower bound to the upper
bound~\eqref{eqn:combined-upper-bound}, there are two factors that are
not perfectly matched.  First, the upper bound depends on $\log
\usedim$, but there is no such dependence in the lower bound. Second,
the upper bound has a term that is proportional to $1/\gamma^2$, but
the corresponding term in the lower bound is proportional
to~$1/\gamma$. Proving a sharper lower bound that closes this gap
remains an open problem.

We remark that Corollary~\ref{coro:lower-bound} follows by a
refinement of the proof of Theorem~\ref{theorem:main-lower-bound}.  In
particular, we first show that the design matrix underlying
Theorem~\ref{theorem:main-lower-bound}---call it
$\NEWXbad$---satisfies the $\gamma_n$-RE condition, where the quantity
$\gamma_n$ converges to zero as a function of sample size~$n$. In
order to prove Corollary~\ref{coro:lower-bound}, we construct a new
block-diagonal design matrix such that each block corresponds to a
version of $\NEWXbad$.  The size of these blocks are then chosen so that,
given a predefined quantity~$\gamma >0$, the new matrix satisfies the
$\gamma$-RE condition.  We then lower bound the prediction error of
this new matrix, using Theorem~\ref{theorem:main-lower-bound} to lower
bound the prediction error of each of the blocks.  We refer the reader
to Section~\ref{sec:proof-coro-lower-bound} for the full proof.

%%%%%%%%%%%%%%%%%%%%%%%%%%%%%%%%%%%%%%%%%%%%%%%%%%%%%%%%%%%%%%%%%%%%%%%%%%

\subsection{Lower bounds for local descent methods}

For any least-squares cost with a coordinate-wise separable
regularizer, Theorem~\ref{theorem:main-lower-bound} establishes the
existence of at least one ``bad'' local minimum such that the
associated prediction error is lower bounded by $1/\sqrt{n}$.  One
might argue that this result could be overly pessimistic, in that the
adversary is given too much power in choosing local minima.  Indeed,
the mere existence of bad local minima need not be a practical concern
unless it can be shown that a typical optimization algorithm will
frequently converge to one of them.

Steepest descent is a standard first-order algorithm for minimizing a
convex cost function~\cite{Bertsekas_nonlin,Boyd04}. However, for
non-convex and non-differentiable loss functions, it is known that the
steepest descent method does not necessarily yield convergence to a 
local minimum~\cite{dem1990introduction,wolfe1975method}.  Although there
exist provably convergent first-order methods for general non-convex
optimization
(e.g.,~\cite{mifflin1982modification,kiwiel1983aggregate}), the paths
defined by their iterations are difficult to characterize, and it is
also difficult to predict the point to which the algorithm eventually
converges.

In order to address a broad class of methods in a unified manner, we
begin by observing that most first-order methods can be seen as
iteratively and approximately solving a local minimization
problem. For example, given a stepsize parameter $\stepsize > 0$, the
method of steepest descent iteratively approximates the minimizer of
the objective over a ball of radius $\stepsize$.  Similarly, the
convergence of algorithms for non-convex optimization algorithms is
based on the fact that they guarantee decrease of the function value
in the local neighborhood of the current
iterate~\cite{mifflin1982modification,kiwiel1983aggregate}.  We thus
study an iterative local descent algorithm taking the form:
\begin{align}
\label{eqn:local-descent-update-formula}
\theta^{t+1} \in \arg \min_{\theta \in \SPECBALL} L(\theta; \lambda),
\end{align}
where $\stepsize > 0$ is a given parameter, and $\Ball_2(\stepsize;
\theta^t) \defn \{ \theta \in \real^\usedim \, \mid \, \ltwos{\theta -
  \theta^t} \leq \stepsize \}$ is the ball of radius $\stepsize$
around the current iterate.  If there are multiple points achieving
the optimum, the algorithm chooses the one that is closest to
$\theta^t$, resolving any remaining ties by randomization.  The
algorithm terminates when there is a minimizer belonging to the
interior of the ball $\SPECBALL$---that is, exactly when
$\theta^{t+1}$ is a local minimum of the loss function.

It should be noted that the
algorithm~\eqref{eqn:local-descent-update-formula} defines a powerful
algorithm---one that might not be easy to implement in polynomial
time---since it is guaranteed to return the global minimum of a
nonconvex program over the ball $\SPECBALL$.  In a certain sense, it
is more powerful than any first-order optimization method, since it
will always decrease the function value at least as much as a descent
step with stepsize related to $\stepsize$.  Since we are proving lower
bounds, these observations only strengthen our result.  We
impose two additional conditions on the regularizers:
\begin{enumerate}[(i)]
  \item[(iv)] Each component function $\reg_j$ is continuous at the
    origin.
  \item[(v)] There is a constant $\curbound$ such that
    $|\reg'_j(x)-\reg'_j(\xtil)|\leq \curbound |x - \xtil|$ for any
    pair $x, \xtil \in (0,\infty)$.
\end{enumerate}
Assumptions (i)-(v) are more restrictive than assumptions (i)-(iii),
but they are satisfied by many popular penalties.  As illustrative
examples, for the $\ell_1$-norm, we have $ H=0$.  For the SCAD
penalty, we have $\curbound = 1/(a-1)$, whereas for the MCP
regularizer, we have $\curbound = 1/b$.  Finally, in order to prevent the
update~\eqref{eqn:local-descent-update-formula} being so powerful that
it reaches the global minimum in one single step, we impose an
additional condition on the stepsize, namely that
\begin{align}
\label{eqn:beta-constraint}
\stepsize \leq \min \Big\{ B, \frac{B}{\lambda \curbound}\Big\}, \quad
\mbox{where $B \defeq \frac{\sigma}{4\sqrt{\numobs}}$.}
\end{align}
It is reasonable to assume that the stepsize bounded by a time-invariant
constant, as we can always partition a single-step update into a finite 
number of smaller steps, increasing the algorithm's time complexity by a
multiplicative constant.  On the other hand, the $\order(1/\sqrt{n})$ 
stepsize is adopted by popular first-order methods.  Under these assumptions, 
we have the following theorem, which applies to any regularizer $\reg$ that 
satisfies Assumptions (i)-(v).

\begin{theorem}
\label{theorem:gradient-descent-lower-bound}
For any pair $(\numobs, \usedim)$ such that $\usedim\geq \numobs \geq
4$, integer $\kdim \geq 2$ and any scalars $\INITVAR \geq 0$ and
$\radius \geq {\sigma}/{\sqrt{\numobs}}$, there is a design matrix
$\Xmat \in \R^{\numobs\times \usedim}$ satisfying the column
normalization condition~\eqref{eqn:column-normalization-condition}
such that
\begin{enumerate}[(a)]
\item The update~\eqref{eqn:local-descent-update-formula} terminates
  after a finite number of steps $\FINIT$ at a vector $\thetahack =
  \theta^{\FINIT+1}$ that is a local minimum of the loss function.
\item Given a random initialization $\theta^0\sim N(0, \INITVAR^2
  I_{d\times d})$, the local minimum satisfies the lower bound
\end{enumerate}
\begin{align*}
\sup_{\thetastar\in \Ball_0(\kdim)\cap \Ball_1(\radius)} \E \left[
  \inf_{\lambda\geq 0}\frac{1}{\numobs} \ltwos{\Xmat\thetahack -
    \Xmat\thetastar}^2 \right] \geq \UNICON \min\{ \radius, \sigma\}
\frac{\sigma}{\sqrt{\numobs}}.
\end{align*}
\end{theorem}

\noindent Theorem~\ref{theorem:gradient-descent-lower-bound} shows
that local descent methods based on a random initialization do not
lead to local optima that achieve the fast rate.  This conclusion
provides stronger negative evidence than
Theorem~\ref{theorem:main-lower-bound}, since it shows that bad local
minima not only exist, but are difficult to avoid.

%%%%%%%%%%%%%%%%%%%%%%%%%%%%%%%%%%%%%%%%%%%%%%%%%%%%%%%%%%%%%%%%%%%%%%%%%

\subsection{Simulations}
\label{SecSimulations}

In the proof of Theorem~\ref{theorem:main-lower-bound} and
Theorem~\ref{theorem:gradient-descent-lower-bound}, we construct
specific design matrices to make the problem hard to solve. In this
section, we apply several popular algorithms to the solution of the sparse linear
regression problem on these ``hard'' examples, and compare their
performance with the $\ell_0$-based
estimator~\eqref{EqnDefnEllZeroEstimator}.  More specifically,
focusing on the special case $\numobs = \usedim$, we perform
simulations for the design matrix $\Xmat \in \R^{\numobs \times
  \numobs}$ used in the proof of
Theorem~\ref{theorem:gradient-descent-lower-bound}.  It is given by
\begin{align*}
\Xmat = \Big[ \mbox{blkdiag} \underbrace{\big \{ \sqrt{\numobs} A,
    \sqrt{\numobs} A, \ldots, \sqrt{\numobs} A \big
    \}}_{\mbox{$\numobs/2$ copies}}\Big],
\end{align*}
where the sub-matrix $A$ takes the form
\begin{align*}
  A = \begin{bmatrix} \cos(\alpha) & -\cos(\alpha)\\ \sin(\alpha) &
    \sin(\alpha)
  \end{bmatrix}, \quad \mbox{where} \quad \alpha = \arcsin(n^{-1/4}).
\end{align*}
Given the $2$-sparse regression vector $\thetastar = \big(0.5, 0.5, 0,
\ldots, 0 \big)$, we form the response vector \mbox{$\yvec = \Xmat
  \thetastar + w$,} where $w \sim N(0, I_{\numobs \times \numobs})$.

We compare the $\ell_0$-based estimator, referred to as the
\emph{baseline estimator}, with three other methods: the Lasso
estimator~\cite{tibshirani1996regression}, the estimator based on the
SCAD penalty~\cite{fan2001variable} and the estimator based on the MCP
penalty~\cite{zhang2010nearly}. In implementing the $\ell_0$-based
estimator, we provide it with the knowledge that $k = 2$, since the
true vector $\thetastar$ is $2$-sparse.  For Lasso, we adopt the
\texttt{MATLAB} implementation~\cite{matlabLasso}, which generates a
Lasso solution path evaluated at $100$ different regularization parameters, 
and we choose the estimate that yields the smallest prediction error.  
For the SCAD penalty, we choose $a = 3.7$ as suggested by Fan and
Li~\cite{fan2001variable}.  For the MCP penalty, we choose $b = 2.7$,
so that the maximum concavity of the MCP penalty matches that of the
SCAD penalty. For the SCAD penalty and the MCP penalty (and recalling
that $\usedim = \numobs$), we studied choices of the regularization
weight of the form $\lambda = C\sqrt{\frac{\log n}{n}}$ for a
pre-factor $C$ to be determined.  As shown in past work on non-convex
regularizers~\cite{loh2013regularized}, such choices of $\lambda$ lead
to low $\ell_2$-error.  By manually tuning the parameter $C$ to
optimize the prediction error, we found that $C=0.1$ is a reasonable
choice.  We used routines from the GIST package~\cite{gong2013gist} to
optimize these non-convex objectives.
\begin{figure}
\centering
\includegraphics[width = 0.7\textwidth]{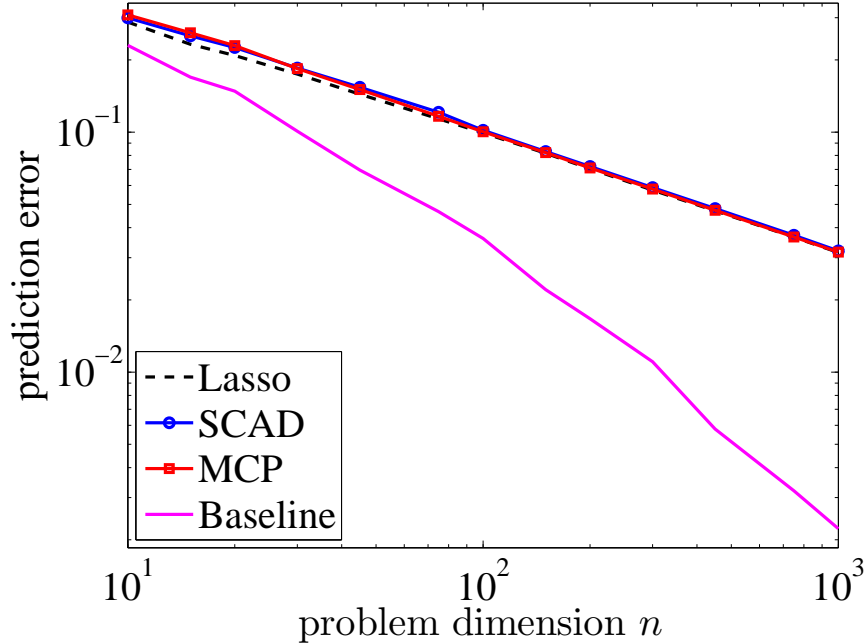}
\caption{Problem scale $n$ versus the prediction error
  $\Exs[\frac{1}{\numobs}\ltwos{\Xmat(\thetahat - \thetastar)}^2]$.
  The expectation is computed by averaging $100$ independent runs of
  the algorithm. Both the sample size $\numobs$ and the prediction
  error are plotted on a logarithmic scale.}
\label{fig:simulation}
\end{figure}

By varying the sample size over the range $10$ to $1000$, we obtained
the results plotted in Figure~\ref{fig:simulation}, in which the
prediction error $\Exs[\frac{1}{\numobs}\ltwos{\Xmat(\thetahat -
    \thetastar)}^2]$ and sample size $\numobs$ are both plotted on a
logarithmic scale. The performance of the Lasso, SCAD-based estimate,
and MCP-based estimate are all similar.  For all of the three methods,
the prediction error scales as $1/\sqrt{n}$, as confirmed by the
slopes of the corresponding lines in Figure~\ref{fig:simulation},
which are very close to~$0.5$.  In fact, by examining the estimator's
output, we find that in many cases, all three estimators output
$\thetahat = 0$, leading to the prediction error
$\frac{1}{n}\ltwos{\Xmat(0 - \thetastar)}^2 = \frac{1}{\sqrt{n}}$.
Since the regularization parameters have been chosen to optimize the
prediction error, this scaling is the best rate that the three
estimators are able to achieve, and it matches the theoretical
prediction of Theorem~\ref{theorem:main-lower-bound} and
Theorem~\ref{theorem:gradient-descent-lower-bound}.

In contrast, the $\ell_0$-based estimator achieves a substantially
better error rate.  The slope of the corresponding line in
Figure~\ref{fig:simulation} is very close to $1$. It means that the
prediction error of the $\ell_0$-based estimator scales as $1/n$,
thereby matching the theoretically-predicted
scaling~\eqref{eqn:l0-optimal-rate}.

%%%%%%%%%%%%%%%%%%%%%%%%%%%%%%%%%%%%%%%%%%%%%%%%%%%%%%%%%%%%%%%%%%%\

\section{Proofs}
\label{SecProofs}

We now turn to the proofs of our theorems and corollary.  In each
case, we defer the proofs of more technical results to the appendices.

%%%%%%%%%%%%%%%%%%%%%%%%%%%%%%%%%%%%%%%%%%%%%%%%%%%%%%%%%%%%%%%%%%%%
 
\subsection{Proof of Theorem~\ref{theorem:main-lower-bound}}
\label{sec:prove-main-result}

For a given triple $(\numobs, \sigma, \radius)$, we define the angle
$\alpha \defeq \arcsin\left(
\frac{\sqrt{\sigma}}{\numobs^{1/4}\sqrt{32 \radius}} \right)$, and the
two-by-two matrix
\begin{subequations}
\begin{align}
\label{eqn:define-matrix-A}
  A = \begin{bmatrix} \cos(\alpha) & -\cos(\alpha)\\ \sin(\alpha) &
    \sin(\alpha)
  \end{bmatrix}.
\end{align}
Using the matrix $A \in \real^{2 \times 2}$ as a building block, we
construct a design matrix $\Xmat \in \real^{\numobs \times
  \usedim}$.  Without loss of generality, we may assume that $\numobs$
is divisible by two. (If $\numobs$ is not divisible by two,
constructing a $(\numobs-1)$-by-$\usedim$ design matrix concatenated
by a row of zeros only changes the result by a constant.) We then
define
\begin{align}
\label{eqn:define-design-matrix-X}
\Xmat = \Big[ \mbox{blkdiag} \underbrace{\big \{ \sqrt{\numobs} A,
    \sqrt{\numobs} A, \ldots, \sqrt{\numobs} A \big
    \}}_{\mbox{$\numobs/2$ copies}}~~{\bf 0}\Big] \in
\real^{\numobs\times \usedim},
\end{align}
\end{subequations}
where the all-zeroes matrix on the right side has dimensions
$\numobs\times (\usedim-\numobs)$.  It is easy to verify that the
matrix $\Xmat$ defined in this way satisfies the column normalization
condition~\eqref{eqn:column-normalization-condition}.

Next, we prove the lower bound~\eqref{eqn:main-lower-bound}.  For any
integers $i,j \in [\usedim]$ with $i < j$, let $\theta_i$ denote the
$i^{th}$ coordinate of $\theta$, and let $\theta_{i:j}$ denote the
subvector with entries $\{\theta_i, \ldots, \theta_j\}$.  Since the
matrix $A$ appears in diagonal blocks of $\Xmat$, we have
\begin{align}
\label{eqn:decompose-prediction-error}
\inf_{\lambda \geq 0} \sup_{\theta\in \Thetahat_\lambda}
\PREDERRSQ{\theta}{\thetastar} & = \inf_{\lambda \geq 0}
\sup_{\theta\in \Thetahat_\lambda} \sum_{i=1}^{\numobs/2}\norms{A
  \big( \theta_{(2i-1):2i} - \opt_{(2i-1):2i} \big)}_2^2
\end{align}
and it suffices to lower bound the right-hand side of the above
equation.

For the sake of simplicity, we introduce the shorthand $B \defeq \frac{4\sigma}{\sqrt{\numobs}}$,
and define the scalars
\begin{align}
\label{eqn:define-gamma-i}
  \gamma_i = \min\{\reg_{2i-1}(B), \reg_{2i}(B)\} \qquad \mbox{for
    each $i = 1, \ldots, \numobs/2$.}
\end{align}
Furthermore, we define
\begin{align}
\label{EqnAdefinition}
a_i \defeq \left\{\begin{array}{ll}
	(\cos \alpha, \sin \alpha) & \mbox{if }\gamma_i = \reg_{2i-1}(B)\\
	(-\cos \alpha, \sin \alpha) & \mbox{if }\gamma_i = \reg_{2i}(B)
\end{array}\right.
\quad \mbox{and} \quad \wprime_i \defeq
\frac{\inprod{a_i}{w_{(2i-1):2i}}}{\sqrt{\numobs}}.
\end{align}
Without loss of generality, we may assume that \mbox{$\gamma_1 = \max
  \limits_{i\in[\numobs/2]}\{\gamma_i\}$} and \mbox{$\gamma_i =
  \reg_{2i -1}(B)$} for all $i \in [\numobs/2]$. If this condition
does not hold, we can simply re-index the columns of $\Xmat$ to make
these properties hold.  Note that when we swap the columns $2i-1$ and
$2i$, the value of $a_i$ doesn't change; it is always associated with
the column whose regularization term is equal to $\gamma_i$.

Finally, we define the regression vector $\thetastar = \begin{bmatrix}
  \frac{\radius}{2} & \frac{\radius}{2} & 0 & \cdots & 0
\end{bmatrix} \in \real^\usedim$.
Given these definitions, the following lemma lower bounds each term on
the right-hand side of
equation~\eqref{eqn:decompose-prediction-error}.

\begin{lemma}
\label{LEMMA:DECOMPOSED-TERM-LOWER-BOUND}
For any $\lambda\geq 0$, there is a local minimum $\thetahat_\lambda$
of the objective function $L(\theta;\lambda)$ such that
$\PREDERRSQ{\thetahatlam}{\thetastar}  \geq \Term_1 + \Term_2$,
where
\begin{subequations}
\begin{align}
\Term_1 & \defeq \indicator\Big[\lambda\gamma_1 > 4B(\sin^2(\alpha)
  \radius \frac{\ltwos{w_{1:2}}}{\sqrt{\numobs}})\Big]
\sin^2(\alpha)(\radius -2B)_+^2 \quad \mbox{and} \\
\Term_2 & \defeq \sum_{i=2}^{\numobs/2} \indicator\Big[B/2 \leq
  \wprime_i \leq B\Big] \Big( \frac{B^2}{4} - \lambda \gamma_1 \Big).
\end{align}
\end{subequations}
Moreover, if the regularizer $\reg$ is convex, then every minimizer
$\thetahat_\lambda$ satisfies this lower bound.
\end{lemma}
\noindent See Appendix~\ref{AppDecomposedOne} for the proof of this
claim.

Using Lemma~\ref{LEMMA:DECOMPOSED-TERM-LOWER-BOUND}, we can now
complete the proof of the theorem.  It is convenient to condition on
the event $\Event \defn \{ \ltwo{w_{1:2}} \leq \frac{\sigma}{32} \}$.
Since $\ltwo{w_{1:2}}^2/\sigma^2$ follows a chi-square distribution
with two degrees of freedom, we have $\mprob[\Event] > 0$.
Conditioned on this event, we now consider two separate cases:

\paragraph{Case 1:}  First, suppose that
$\lambda\gamma_1 > \sigma^2/\numobs$.  In this case, we have
\begin{align*}
  4 B \Big \{ \sin^2(\alpha) \radius +
  \frac{\ltwos{w_{1:2}}}{\sqrt{\numobs}} \Big \} & \leq
  \frac{16\sigma}{\sqrt{\numobs}} \left(
  \frac{\sigma}{32\sqrt{\numobs}} + \frac{\sigma}{32\sqrt{\numobs}}
  \right) = \frac{\sigma^2}{\numobs} < \lambda\gamma_1,
\end{align*}
and consequently
\begin{subequations}
\begin{align}
\label{EqnCaseOneBound}
\Term_1 + \Term_2 \geq \Term_1 = \sin^2(\alpha)(\radius -2B)_+^2 =
\frac{\sigma}{32\sqrt{\numobs}\radius} \left(\radius -
\frac{4\sigma}{\sqrt{\numobs}}\right)^2 \geq \frac{\sigma
  \radius}{128\sqrt{\numobs}},
\end{align}
where the last inequality holds since we have assumed that $\radius
\geq 8 \sigma / \sqrt{\numobs}$.

\paragraph{Case 2:} 
Otherwise, we may assume that $\lambda\gamma_1 \leq \sigma^2/\numobs$.
In this case, we have
\begin{align}
\label{EqnCaseTwoBound}
\Term_1 + \Term_2 & \geq \Term_2 = \sum_{i=2}^{\numobs/2}
\indicator\left(B/2 \leq \wprime_i \leq B\right)
\frac{3\sigma^2}{\numobs}.
\end{align}
\end{subequations}

\vspace*{.1in}

Combining the two lower bounds~\eqref{EqnCaseOneBound}
and~\eqref{EqnCaseTwoBound}, we find
\begin{align*}
&\E\left[ \inf_{\lambda \geq 0} \sup_{\theta\in \Thetahat_\lambda}
  \frac{1}{\numobs} \ltwos{\Xmat\theta - \Xmat\thetastar}^2 \Big|
  \event \right] \\
  & \qquad \geq \underbrace{\E\left[ \min\left\{ \frac{\sigma
      \radius}{128 \sqrt{\numobs}}, \sum_{i=2}^{\numobs/2}
    \indicator\left[B/2 \leq \wprime_i \leq B/2\right]
    \frac{3\sigma^2}{\numobs} \right\}\right]}_{\Term_3},
\end{align*}
where we have used the fact that $\{\wprime_i\}_{i=2}^{\numobs/2}$ are
independent of the event $\Event$.  Using the inequality
$\min\left\{a, \sum_{i=2}^{n/2} b_i\right\} \geq \sum_{i=2}^{n/2}
\min\{2 a/n, b_i\}$, valid for scalars $a$ and
$\{b_i\}_{i=2}^{\numobs/2}$, we see that
\begin{align*}
\Term_3 & \geq \sum_{i=2}^{\numobs/2} \mprob\left[ \frac{2
    \sigma}{\sqrt{\numobs} } \leq \wprime_i \leq \frac{4
    \sigma}{\sqrt{\numobs}} \right] \min\left\{ \frac{\sigma
  \radius}{128 \numobs \sqrt{\numobs}}, \frac{3\sigma^2}{\numobs}
\right\},
\end{align*}
where we have used the fact that $\Exs\Big[\indicator \big[B/2 \leq
    \wprime_i \leq B \big] \Big] = \mprob[B/2 \leq \wprime_i \leq B]$,
and the definition $B \defeq \frac{4 \sigma}{\sqrt{\numobs}}$.

Since $\wprime_i \sim N(0,\sigma^2/\numobs)$, the probability $\mprob[
  {2\sigma}/{\sqrt{\numobs}} \leq \wprime_i \leq
  {4\sigma}/{\sqrt{\numobs}} ]$ is bounded away from zero
independently of all problem parameters. Hence, there is a universal
constant $c_2 > 0$ such that $\Term_3 \geq c_2 \min\left\{
\frac{\sigma \radius}{\sqrt{\numobs}}, \; \sigma^2 \right\}$.  Putting
together the pieces, we have shown that
\begin{align*}
\E\left[ \inf_{\lambda \geq 0} \sup_{\theta\in \Thetahat_\lambda}
  \frac{1}{\numobs} \ltwos{\Xmat\theta - \Xmat\thetastar}^2 \right]
&\geq \: \mprob[\event] \, \Term_3 \; \geq c_2' \min \left \{ \frac{
  \sigma \radius}{\sqrt{\numobs}}, \: \sigma^2 \right\},
\end{align*}
which completes the proof of the theorem.

%%%%%%%%%%%%%%%%%%%%%%%%%%%%%%%%%%%%%%%%%%%%%%%%%%%%%%%%%%%%%%%%%%%%%%

\subsection{Proof of Corollary~\ref{coro:lower-bound}}
\label{sec:proof-coro-lower-bound}

Here we provide a detailed proof of
inequality~\eqref{eqn:coro-lower-bound}.  We note that
inequality~\eqref{eqn:coro-lower-bound-for-convex-penalty} follows by
an essentially identical series of steps, so that we omit the details.

Let $m$ be an even integer and let $X_m \in \R^{m\times m}$ denote the
design matrix constructed in the proof of
Theorem~\ref{theorem:main-lower-bound}. In order to avoid confusion,
we rename the parameters $(n,d,R)$ in the
construction~\eqref{eqn:define-design-matrix-X} by $(n',d',R')$, and
set them equal to
\begin{align}
\label{eqn:m-m-radius}
(n',d',R') \defeq \Big(m, m, \min \Big \{ \frac{R \sqrt{n}}{k
  \sqrt{m}}, \frac{\sigma}{16 \gamma \sqrt{m}} \Big\} \Big),
\end{align}
where the quantities $(k,m,n,R,\sigma)$ are defined in the statement
of Corollary~\ref{coro:lower-bound}.  Note that $X_m$ is a square
matrix, and according to equation~\eqref{eqn:define-design-matrix-X},
all of its eigenvalues are lower bounded by
$(\frac{m^{1/2}\sigma}{16R'})^{1/2}$. By
equation~\eqref{eqn:m-m-radius}, this quantity is lower bounded by
$\sqrt{m\gamma}$.

Using the matrix $X_m$ as a building block, we now construct a larger
design matrix \mbox{$X \in \R^{n \times n}$} that we then use to prove
the corollary. Let $m$ be the greatest integer divisible by two such
that $k m \leq n$. By the assumption that $n \geq 4 k$, we have $m
\geq 4$. Consequently, we may construct the $n \times n$ dimensional
matrix
\begin{align}
\label{eqn:define-coro-design-matrix-X}
\Xmat \defeq \mbox{blkdiag} \Big \{\underbrace{ \sqrt{\numobs/m} X_m,
  \ldots, \sqrt{\numobs/m} X_m}_{\mbox{$k$ copies}}, \sqrt{n} I_{n-km}
\Big\} \in \real^{\numobs\times \numobs},
\end{align}
where $I_{n-km}$ is the $(n-km)$-dimensional identity matrix. It is
easy to verify the matrix $X$ satisfies the column normalization
condition. Since all eigenvalues of $X_m$ are lower bounded by
$\sqrt{m\gamma}$, we are guaranteed that all eigenvalues of $X$ are
lower bounded by $\sqrt{n\gamma}$. Thus, the matrix $X$ satisfies the
$\gamma$-RE condition.

It remains to prove a lower bound on the prediction error, and in
order to do so, it is helpful to introduce some shorthand
notation. Given an arbitrary vector $u \in \R^n$, for each integer $i
\in \{1, \ldots, k\}$, we let $u_{(i)} \in \real^m$ denote the
sub-vector consisting of the $((i-1)m+1)$-th to the $(im)$-th elements
of vector $u$, and we let $u_{(k+1)}$ denote the sub-vector consisting
of the last $n-km$ elements. We also introduce similar notation for
the function $\rho(x) = \rho_1(x_1) + \dots +\rho_n(x_n)$;
specifically, for each $i \in \{1, \ldots, k\}$, we define the
function $\rho_{(i)}: \R^m\to \R$ via $\rho_{(i)}(\theta) \defeq
\sum_{j=1}^m \rho_{(i-1)m+j}(\theta_j)$.

Using this notation, we may rewrite the cost function as:
\begin{align*}
L(\theta;\lambda) = \frac{1}{n}\sum_{i=1}^k\Big(\ltwos{\sqrt{n/m}X_m
  \theta_{(i)} - y_{(i)}}^2 + n \lambda \rho_{(i)}(\theta_{(i)}) \Big)
+ h(\theta_{(k+1)}),
\end{align*}
where $h$ is a function that only depends on $\theta_{(k+1)}$.  If we
define $\theta'_{(i)} \defeq \sqrt{n/m}\theta_{(i)}$ and
$\rho'_{(i)}(\theta) \defeq \frac{n}{m}\rho_{(i)}(\sqrt{m/n}\theta)$,
then substituting them into the above expression, the cost function
can be rewritten as
\begin{align*}
G(\theta'; \lambda) & \defeq \frac{m}{n} \sum_{i=1}^k\Big(\frac{1}{m}
\ltwos{X_m \theta'_{(i)} - y_{(i)}}^2 + \lambda
\rho'_{(i)}(\theta'_{(i)}) \Big) + h(\sqrt{m/n} \; \theta'_{(k+1)}).
\end{align*}
Note that if the vector $\thetahat$ is a local minimum of the function
$\theta \mapsto L(\theta;\lambda)$, then the rescaled vector
$\thetahat' \defeq \sqrt{n/m}\;\thetahat$ is a local minimum of the
function $\theta' \mapsto G(\theta';\lambda)$.  Consequently, the
sub-vector $\thetahat'_{(i)}$ must be a local minimum of the function
\begin{align}
\label{eqn:coro-split-local-minimum}
\frac{1}{m}\ltwos{X_m \theta'_{(i)} - y_{(i)}}^2 +
\rho'_{(i)}(\theta'_{(i)}).
\end{align}
Thus, the sub-vector $\thetahat'_{(i)}$ is the solution of a
regularized sparse linear regression problem with design matrix
$X_m$. 

Defining the rescaled true regression vector $(\thetastar)' \defeq
\sqrt{n/m}\;\thetastar$, we can then write the prediction error as
\begin{align}
\frac{1}{n} \ltwos{X(\thetahat - \thetastar)}^2 & =
\frac{1}{n}\sum_{i=1}^k \Big(\ltwos{X_m(\thetahat'_{(i)} -
  (\thetastar)'_{(i)})}^2 \Big) + \ltwos{\thetahat_{(k+1)} -
  \thetastar_{(k+1)}}^2 \nonumber \\
\label{eqn:coro-split-prediction-error}
& \geq \frac{m}{n}\sum_{i=1}^k \Big( \frac{1}{m}
\ltwos{X_m(\thetahat'_{(i)} - (\thetastar)'_{(i)})}^2 \Big).
\end{align}
Consequently, the overall prediction error is lower bounded by a
scaled sum of the prediction errors associated with the design
matrix~$X_m$. Moreover, each term $\frac{1}{m}
\ltwos{X_m(\thetahat'_{(i)} - (\thetastar)'_{(i)})}^2$ can be bounded
by Theorem~\ref{theorem:main-lower-bound}.

More precisely, let $\mathcal{Q}(X,2k,R)$ denote the left-hand side of
inequality~\eqref{eqn:coro-lower-bound}. The above analysis shows that
the sparse linear regression problem on the design matrix $X$ and the
constraint $\thetastar\in \Ball_0(2\kdim)\cap \Ball_1(\radius)$ can be
decomposed into smaller-scale problems on the design matrix $X_m$ and
constraints on the scaled vector $(\thetastar)'$.  By the rescaled
definition of $(\thetastar)'$, the constraint $\thetastar\in
\Ball_0(2\kdim)\cap \Ball_1(\radius)$ holds if and only if
$(\thetastar)' \in \Ball_0(2\kdim)\cap \Ball_1(\sqrt{n/m}\radius)$.
Recalling the definition of the radius $R'$ from
equation~\eqref{eqn:m-m-radius}, we can ensure that $(\thetastar)' \in
\Ball_0(2\kdim)\cap \Ball_1(\sqrt{n/m}\radius)$ by requiring that
$(\thetastar)'_{(i)}\in \Ball_0(2)\cap \Ball_1(R')$ for each
\mbox{index $i \in \{1,\dots, k\}$.}  Combining
expressions~\eqref{eqn:coro-split-local-minimum}
and~\eqref{eqn:coro-split-prediction-error}, the quantity
$\mathcal{Q}(X,2k,R)$ can be lower bounded by the sum
\begin{subequations}
\begin{align}
\label{eqn:bound-q-X2kR}
\mathcal{Q}(X,2k,R) \geq \frac{m}{n}\sum_{i=1}^k
\mathcal{Q}(X_m,2,R').
\end{align}
By Theorem~\ref{theorem:main-lower-bound}, we have
\begin{align}
\label{eqn:bound-q-Xm2Rp}
\mathcal{Q}(X_m,2,R') \geq c \min \left\{ \sigma^2,
\frac{\sigma\radius'}{\sqrt{m}}\right\} = c \min \left\{\sigma^2,
\frac{\sigma^2}{16\gamma m}, \frac{\sigma \radius \sqrt{n}}{k m}
\right\},
\end{align}
\end{subequations}
where the second equality follows from our choce of $R'$ from
equation~\eqref{eqn:m-m-radius}.  Combining the lower
bounds~\eqref{eqn:bound-q-X2kR} and~\eqref{eqn:bound-q-Xm2Rp}
completes the proof.

%%%%%%%%%%%%%%%%%%%%%%%%%%%%%%%%%%%%%%%%%%%%%%%%%%%%%%%%%%%%%%%%%%%%%%%%%%%%%%%%

\subsection{Proof of Theorem~2}

The proof of Theorem~2 is conceptually similar to the proof of
Theorem~\ref{theorem:main-lower-bound}, but differs in some key
details.  We begin with the altered definitions
\begin{align*}
  \alpha \defeq \arcsin\left( \frac{\sqrt{\sigma}}{\numobs^{1/4}
    \sqrt{\altradius}} \right) \quad \mbox{and} \quad B \defeq
  \frac{\sigma}{4\sqrt{\numobs}}, \qquad \mbox{where } \altradius
  \defeq \min\{\radius,\sigma\}.
\end{align*}
Given our assumption $R\geq \sigma/\sqrt{\numobs}$, note that we are
guaranteed that the inequality \mbox{$2B = \sigma/(2\sqrt{\numobs})
  \leq r/2$} holds.  We then define the matrix $A\in \R^{2\times2}$
and the matrix $\Xmat\in \R^{\numobs\times \usedim}$ by
equations~\eqref{eqn:define-matrix-A}
and~\eqref{eqn:define-design-matrix-X}.

%%%%%%%%%%%%%%%%%%%%%%%%%%%%%%%%%%%%%%%%%%%%%%%%%%%%%%%%%%%%%%%%%%%%%%%%%

\subsubsection{Proof of part (a)}

Let $\{\theta^t\}_{t=0}^\infty$ be the sequence of iterates generated
by equation~\eqref{eqn:local-descent-update-formula}.  We proceed via
proof by contradiction, assuming that the sequence does not terminate
finitely, and then deriving a contradiction.   We begin with a lemma. 

\begin{lemma} 
\label{LemAuxiliary}
If the sequence of iterates $\{\theta^t\}_{t=0}^\infty$ is not
finitely convergent, then it is unbounded.
\end{lemma}

We defer the proof of this claim to the end of this section.  Based on
Lemma~\ref{LemAuxiliary}, it suffices to show that, in fact, the
sequence $\{\theta^t\}_{t=0}^\infty$ is bounded.  Partitioning the
full vector as \mbox{$\theta \defeq \big(\theta_{1:n}, \theta_{n+1:d}
  \big)$,} we control the two sequences $\{\theta^t_{1:n}
\}_{t=0}^\infty$ and $\{ \theta^t_{n+1:d} \}_{t=0}^\infty$. \\

Beginning with the former sequence, notice that the objective function
can be written in the form
\begin{align*}
L(\theta;\lambda) = \frac{1}{\numobs} \ltwos{y -
  X_{1:\numobs}\theta_{1:n}}^2 + \sum_{i=1}^\usedim \lambda
\rho_i(\theta_i),
\end{align*}
where $\Xmat_{1:\numobs}$ represents the first $\numobs$ columns of
matrix $\Xmat$. The conditions~\eqref{eqn:define-matrix-A}
and~\eqref{eqn:define-design-matrix-X} guarantee that the Gram matrix
$X_{1:n}^T X_{1:n}$ is positive definite, which implies that the
quadratic function \mbox{$\theta_{1:n} \mapsto \ltwos{y -
    X_{1:n}\theta_{1:n}}^2$} is strongly convex.  Thus, if the
sequence $\{\theta_{1:n}^t\}_{t=0}^\infty$ were unbounded, then the
associated cost sequence $\{L(\theta^t; \lambda)\}_{t=0}^\infty$ would
also be unbounded.  But this is not possible since $L(\theta^t;
\lambda) \leq L(\theta^0; \lambda)$ for all iterations $t = 1, 2,
\ldots$.  Consequently, we are guaranteed that the sequence
$\{\theta^t_{1:n}\}_{t=0}^\infty$ must be bounded.

It remains to control the sequence
$\{\theta_{n+1:d}^t\}_{t=0}^\infty$.  We claim that for any $i \in
\{n+1, \ldots, \usedim\}$, the sequence
$\{|\theta_i^t|\}_{t=0}^\infty$ is non-increasing, which implies the
boundedness condition.  Proceeding via proof by contradiction, suppose
that $|\theta^t_i| < |\theta^{t+1}_i|$ for some index $i \in \{n+1,
\ldots, \usedim \}$ and iteration number $t \geq 0$.  Under this
condition, define the vector
\begin{align*}
	\widetilde \theta^{t+1}_j \defeq \begin{cases} \theta^{t+1}_j
          & \mbox{if $j\neq i$}\\ \theta^{t}_j & \mbox{if $j=
            i$.}  \end{cases}
\end{align*}
Since $\rho_j$ is a monotonically non-decreasing function of $|x|$, we
are guaranteed that $L(\widetilde \theta^{t+1};\lambda) \leq
L(\theta^{t+1};\lambda)$, which implies that $\widetilde \theta^{t+1}$
is also a constrained minimum point over the ball $\SPECBALL$.  In
addition, we have
\begin{align*}
\ltwos{\widetilde \theta^{t+1} - \theta^{t}} = \ltwos{\theta^{t+1} -
  \theta^{t}} - |\theta^t_i - \theta^{t+1}_i| < \stepsize,
\end{align*}
so that $\widetilde \theta^{t+1}_j$ is strictly closer to $\theta^t$.
This contradicts the specification of the algorithm, in that it
chooses the minimum closest to $\theta^t$.

\paragraph{Proof of Lemma~\ref{LemAuxiliary}:}

The final remaining step is to prove Lemma~\ref{LemAuxiliary}.  We
first claim that $\ltwos{\theta^s - \theta^t} \geq \stepsize$ for all
pairs $s < t$.  If not, we could find some pair $s < t$ such that
$\ltwos{\theta^s - \theta^t} < \stepsize$.  But since $t > s$, we are
guaranteed that $L(\theta^t;\lambda) \leq L(\theta^{s+1};\lambda)$.
Since $\theta^{s+1}$ is a global minimum over the ball
$\Ball_2(\stepsize; \theta^s)$ and $\ltwos{\theta^s - \theta^t} <
\stepsize$, the point $\theta^t$ is also a global minimum, and this
contradicts the definition of the algorithm (since it always chooses
the constrained global minimum closest to the current iterate).

Using this property, we now show that $\{\theta^t\}_{t=0}^\infty$ is
unbounded.  For each \mbox{iteration} \mbox{$t = 0, 1, 2 \ldots$,} we use
$\Ball^t = \Ball_2(\stepsize/3; \theta^t)$ to denote the Euclidean
ball of radius $\stepsize/3$ centered at $\theta^t$.  Since
$\ltwos{\theta^s - \theta^t} \geq \stepsize$ for all $s \neq t$, the
balls $\{\Ball^t\}_{t=0}^\infty$ are all disjoint, and hence there is
a numerical constant $C > 0$ such that for each $T \geq 1$, we have
\begin{align*}
{\rm vol}\Big( \cup_{t=0}^T \Ball^t \Big) = \sum_{t=0}^T {\rm
  vol}(\Ball^t) = C \sum_{t=0}^T\stepsize^d.
\end{align*}
Since this volume diverges as $T \rightarrow \infty$, we conclude that
the set $\Ball \defn \cup_{t=0}^\infty \Ball^t$ must be unbounded.  By
construction, any point in $\Ball$ is within $\stepsize/3$ of some
element of the sequence $\{\theta^t\}_{t=0}^\infty$, so this sequence
must be unbounded, as claimed.

%%%%%%%%%%%%%%%%%%%%%%%%%%%%%%%%%%%%%%%%%%%%%%%%%%%%%%%%%%%%%%%%%%%%%%%

\subsubsection{Proof of part (b)}

We now prove a lower bound on the prediction error corresponding the
local minimum to which the algorithm converges, as claimed in part (b)
of the theorem statement.  In order to do so, we begin by introducing
the shorthand notation
\begin{align}
\label{eqn:define-gamma-i}
\gamma_i = \min\{\sup_{u\in (0,B]}\reg'_{2i-1}(u), \sup_{u\in (0,B]}
    \reg'_{2i}(u)\} \qquad \mbox{for each $i = 1, \ldots, \numobs/2$}.
\end{align}
Then we define the quantities $a_i$ and $\wprime_i$ by
equations~\eqref{EqnAdefinition}.  Similar to the proof of
Theorem~\ref{theorem:main-lower-bound}, we assume (without loss of
generality, re-indexing as needed) that $\gamma_i = \sup_{u\in
  (0,B]}\reg'_{2i-1}(u)$ and that $\gamma_1 = \max_{i\in[\numobs/2]}
  \{\gamma_i\}$.

Consider the regression vector $\thetastar \defn \begin{bmatrix}
  \frac{r}{2} & \frac{r}{2} & 0 & \cdots & 0
\end{bmatrix}$.
Since the matrix $A$ appears in diagonal blocks of $\Xmat$, the
algorithm's output $\thetahack$ has error
\begin{align}
\label{eqn:gradient-decomposable-m-estimator}
\inf_{\lambda\geq 0} \PREDERRSQ{\thetahack}{\thetastar} & =
\inf_{\lambda\geq 0} \sum_{i=1}^{\numobs/2} \norms{A
  \big(\thetahack_{(2i-1):2i} - \opt_{(2i-1):2i} \big)}_2^2.
\end{align}
Given the random initialization $\thetait{0}$, we define the events
\begin{align*}
\Event_0 \defeq \Big \{ \max\{ \theta_{1}^0, \theta_{2}^0 \} \leq 0\Big\}
\quad \mbox{and} \quad \event_1 \defeq \Big\{
\lambda\gamma_1 \geq 2\sin^2(\alpha)\altradius +
\frac{2 \ltwos{w_{1:2}}}{\sqrt{\numobs}} + 3B \Big \},
\end{align*}
as well as the (random) subsets
\begin{align*}
\Sset_1 & \defeq \Big \{ i \in \big \{ 2, \ldots, \numobs/2 \big \} \,
\mid \, \lambda \gamma_1 \leq 2\wprime_i - 4 B \Big \}, \quad
\mbox{and} \\
\Sset_2 & \defeq \Big\{ i \in \{2, \ldots, \numobs/2 \} \, \mid \,
\mbox{$2\sin^2(\alpha)\altradius + \frac{2
    \ltwos{w_{1:2}}}{\sqrt{\numobs}} + 3B \leq 2\wprime_i - 4 B$}
\Big\}.
\end{align*}
Here the reader should recall the definition of $\wprime$ from
equation~\eqref{EqnAdefinition}.

Given these definitions, the following lemma provides lower bounds on
the decomposition~\eqref{eqn:gradient-decomposable-m-estimator} for
the vector $\thetahack$ after convergence.

\begin{lemma}
\label{LEMMA:GRADIENT-DESCENT-TERMWISE-LOWER-BOUND}
\begin{enumerate}[(a)]
\item If $\HACKETWO$ holds, then $\norms{A \, \big(\thetahack_{1:2} -
  \opt_{1:2}\big)}_2^2 \geq \frac{\sigma \altradius}{4\sqrt{n}}$.
\item 
For any index $i \in \Sset_1$, we have
$\norms{A \, \big(\thetahack_{2i-1:2i} - \opt_{2i-1:2i} \big)}_2^2 \geq
\frac{\sigma \altradius}{8\numobs^{3/2}}$.
\item We have $\mprob[\event_0] = 1/4$, and moreover $\min
  \limits_{i \in \{2, \ldots, \numobs/2 \}} \mprob[i\in \SsetThree]
  \geq c$ for some numerical constant $c > 0$.
\end{enumerate}
\end{lemma}
\noindent See Appendix~\ref{AppDecomposedTwo} for the proof of this
claim. \\

Conditioned on event $\event_0$, for any index $i \in \Sset_2$,
either the event $\HACKETWO$ holds, or we have
\begin{align*}
\lambda \gamma_1 < 2\sin^2(\alpha)\altradius + \frac{2
  \ltwos{w_{1:2}}}{\sqrt{\numobs}} + 3 B \leq 2\wprime_i - 4 B,
\end{align*}
which means that $i\in
\Sset_1$ holds. Applying
Lemma~\ref{LEMMA:GRADIENT-DESCENT-TERMWISE-LOWER-BOUND} yields the
lower bound
\begin{align*}
  \inf_{\lambda\geq 0} \sum_{i=1}^{\numobs/2} \norms{A
    \thetahack_{(2i-1):2i} - A \opt_{(2i-1):2i}}_2^2 &\geq
  \indicator[\Event_0] \; \min\left\{ \frac{\sigma \altradius}{4
    \sqrt{n}}, \: \; \frac{\sigma \altradius}{8\numobs^{3/2}} \;
  \sum_{i=2}^{[\numobs/2]} \indicator[i \in \SsetThree] \right\}\\ &=
  \indicator[\Event_0]\; \frac{\sigma \altradius}{8\numobs^{3/2}} \;
  \sum_{i=2}^{[\numobs/2]} \indicator[i \in \SsetThree],
\end{align*}
where the last equality holds since $\frac{\sigma
  \altradius}{8\numobs^{3/2}} \; \sum_{i=2}^{[\numobs/2]} \indicator[i
  \in \Sset_2] \leq \frac{\sigma \altradius}{8\numobs^{3/2}}\;(n/2-1)
< \frac{\sigma \altradius}{4 \sqrt{n}}$.  Since the event $\event_0$
is independent of the event $\{i \in \SsetThree, i = 2, \ldots,
\numobs/2 \}$, we have
\begin{align*}
\Exs \Big[ \inf_{\lambda\geq 0} \sum_{i=1}^{\numobs/2} \norms{A
    \thetahack_{(2i-1):2i} - A \opt_{(2i-1):2i}}_2^2 \Big] & \geq
\mprob[\event_0] \; \frac{\sigma \altradius}{8\numobs^{3/2}}
\sum_{i=2}^{[\numobs/2]} \mprob[i \in \SsetThree] \\
& \stackrel{(i)}{\geq} \frac{1}{4} \; \frac{\sigma
  \altradius}{8\numobs^{3/2}} \;c(n/2-1)\\
& = c' \min \{\radius, \sigma \} \, \frac{\sigma}{\sqrt{\numobs}},
\end{align*}
where step (i) uses the lower bound $\mprob[\event_0] = 1/4$ and
$\mprob[i\in \SsetThree] \geq c$ from
Lemma~\ref{LEMMA:GRADIENT-DESCENT-TERMWISE-LOWER-BOUND}.  Combined
with the decomposition~\eqref{eqn:gradient-decomposable-m-estimator},
the proof is complete.

%%%%%%%%%%%%%%%%%%%%%%%%%%%%%%%%%%%%%%%%%%%%%%%%%%%%%%%%%%%%%%%%%%%%%%%%%%%%%%

\section{Discussion}
\label{SecDiscussion}

In this paper, we have demonstrated a fundamental gap in sparse linear
regression: the best prediction risk achieved by a class of
$M$-estimators based on coordinate-wise separable regularizers is
strictly larger than the the classical minimax prediction risk,
achieved for instance by minimization over the $\ell_0$-ball. This gap
applies to a range of methods used in practice, including the Lasso in
its ordinary and weighted forms, as well as estimators based on
nonconvex penalties such as the MCP and SCAD penalties.

Several open questions remain, and we discuss a few of them here.
When the penalty function $\rho$ is convex, the M-estimator minimizing
function~\eqref{EqnLoss} can be understood as a particular convex
relaxation of the $\ell_0$-based
estimator~\eqref{EqnDefnEllZeroEstimator}.  It would be interesting to
consider other forms of convex relaxations for the $\ell_0$-based
problem.  For instance, Pilanci et al.~\cite{PilWaiElg15} show how a
broad class of $\ell_0$-regularized problems can be reformulated
exactly as optimization problems involving convex functions in Boolean
variables.  This exact reformulation allows for the direct application
of many standard hierarchies for Boolean polynomial programming,
including the Lasserre hierarchy~\cite{lasserre2001explicit} as well
as the Sherali-Adams hierarchy~\cite{Sherali90}.  Other relaxations
are possible, including those that are based on introducing auxiliary
variables for the pairwise interactions (e.g., $\gamma_{ij} = \theta_i
\theta_j$), and so incorporating these constraints as polynomials in
the constraint set.  We conjecture that for any fixed natural number
$t$, if the the $t$-th level Lasserre (or Sherali-Adams) relaxation is
applied to such a reformulation, it still does not yield an estimator
that achieves the fast rate~\eqref{eqn:l0-optimal-rate}.  Since a
$t^{th}$-level relaxation involves $\order(\usedim^t)$ variables, this
would imply that these hierarchies do not contain polynomial-time
algorithms that achieve the classical minimax risk.  Proving or
disproving this conjecture remains an open problem.

Finally, when the penalty function $\rho$ is concave, concurrent work
by Ge et al.~\cite{ge2015strong} shows that finding the global minimum
of the loss function~\eqref{EqnLoss} is strongly NP-hard. This result
implies that no polynomial-time algorithm computes the global minimum
unless ${\bf NP}={\bf P}$. The result given here is complementary in
nature: it shows that bad local minima exist, and that local descent
methods converge to these bad local minima. It would be interesting to
extend this algorithmic lower bound to a broader class of first-order
methods. For instance, we suspect that any algorithm that relies on an
oracle giving first-order information will inevitably converge to a
bad local minimum for a broad class of random initializations.

%%%%%%%%%%%%%%%%%%%%%%%%%%%%%%%%%%%%%%%%%%%%%%%%%%%%%%%%%%%%%%%%%%%%%%%%%%%%%%

\subsection*{Acknowledgements}

This work was partially supported by grants NSF grant DMS-1107000, NSF
grant CIF-31712-23800, Air Force Office of Scientific Research Grant
AFOSR-FA9550-14-1-0016, and Office of Naval Research MURI
N00014-11-1-0688.

%%%%%%%%%%%%%%%%%%%%%%%%%%%%%%%%%%%%%%%%%%%%%%%%%%%%%%%%%%%%%%%%%%%%%%%%%%%%%
\appendix

%%%%%%%%%%%%%%%%%%%%%%%%%%%%%%%%%%%%%%%%%%%%%%%%%%%%%%%%%%%%%%%%%%%%%%%%%%%%%%

\section{Fast rate for the bad example of Dalalyan et al.~\cite{dalalyan2014prediction}}
\label{sec:fast-rate-dalalyan}

In this appendix, we describe the bad example of Dalayan et
al.~\cite{dalalyan2014prediction}, and show that a reweighted form of
the Lasso achieves the fast rate.  For a given sample size $\numobs \geq 4$,
they consider a linear regression model $y = X \thetastar + w$, where
$X \in \real^{\numobs \times 2 m}$ with $m = n - 1$, and the noise
vector $w \in \{-1, 1\}^n$ has i.i.d. Rademacher entries (equiprobably
chosen in $\{-1,1\}$).  In the construction, the true vector
$\thetastar$ is 2-sparse, and the design matrix $X\in \R^{n\times 2m}$
is given by
\begin{align*}
	X = \sqrt{n} \begin{bmatrix} {\bf 1}^T_m & {\bf
            1}^T_m\\ I_{m\times m} & - I_{m\times m}
        \end{bmatrix}
\end{align*}
where ${\bf 1}_m \in \real^m$ is a vector of all ones. Notice that
this construction has $n = m+1$.\\

In this appendix, we analyze the performance of the following
estimator
\begin{align}
\label{eqn:loss}
\thetahat \in \arg\min_{\theta\in \R^{2m}} \frac{1}{n}\ltwos{X\theta -
  y}^2 + \lambda \sum_{i=2}^m (|\theta_i| + |\theta_{m+i}|).
\end{align}
It is a reweighted form of the Lasso based on $\ell_1$-norm
regularization, but one that imposes \emph{no} constraint on the first
and the $(m+1)$-th coordinate.  We claim that with an appropriate
choice of $\lambda$, this estimator achieves the fast rate for any
$2$-sparse vector $\thetastar$.

Letting $\thetahat$ be a minimizer of function~\eqref{eqn:loss}, we
first observe that no matter what value it attains, the minimizer
always chooses $\thetahat_1$ and $\thetahat_{m+1}$ so that
$(X\thetahat)_{1:2} = y_{1:2}$. This property occurs because:
\begin{itemize}
\item There is no penalty term associated with $\thetahat_1$ and
  $\thetahat_{m+1}$.
\item By the definition of $X$, changes in the coordinates
  $\thetahat_1$ and $\thetahat_{m+1}$ only affect the first two
  coordinates of $X\thetahat$ by the additive term
\begin{align*}
	\sqrt{n} \begin{bmatrix} 1 & 1\\ e_1 & -e_1
        \end{bmatrix}
\cdot
\begin{bmatrix}
 \thetahat_1\\ \thetahat_{m+1}
\end{bmatrix}
\end{align*}
Since the above 2-by-2 matrix is non-singular, there is always an
assignment to $(\thetahat_1, \thetahat_{m+1})$ so that
$(X\thetahat)_{1:2} - y_{1:2} = 0$.
\end{itemize}
Thus, only the last $n-2$ coordinates of $X\thetahat - y$ might be
non-zero, so that we may rewrite the objective
function~\eqref{eqn:loss} as
\begin{align}\label{eqn:simplified-loss}
	&\frac{1}{n}\ltwos{X\thetahat - y}^2 + \lambda \sum_{i=2}^m
  (|\thetahat_i| + |\thetahat_{m+i}|)\nonumber \\ &\qquad =
  \frac{1}{n} \sum_{i=2}^m \left( (\sqrt{n}\thetahat_i -
  \sqrt{n}\thetahat_{m+i} - y_i)^2 + \lambda (|\thetahat_i| +
  |\thetahat_{m+i}|) \right)
\end{align}
The function~\eqref{eqn:simplified-loss} is not strictly convex so
that there are multiple equivalent solutions.  Essentially, we need to
break symmetry by choosing to vary one of $\thetahat_i$ or
$\thetahat_{m+i}$, for each $i \in \{2, \ldots, m\}$. Without loss of
generality, we assume that $\thetahat_{m+2} = \thetahat_{m+3} = \dots
= \thetahat_{2m} = 0$, so that the equation is simplified as
\begin{align}
\label{eqn:simplified-loss-2}
&\frac{1}{n}\ltwos{X\thetahat - y}^2 + \lambda \sum_{i=2}^m
(|\thetahat_i| + |\thetahat_{m+i}|) = \frac{1}{n} \sum_{i=2}^m \left(
(\sqrt{n}\thetahat_i - y_i)^2 + \lambda |\thetahat_i|) \right).
\end{align}
Moreover, with this choice, we can write the prediction error as
\begin{align}
\label{eqn:simplified-pred-error}
R(\thetahat) \defeq \frac{1}{n}\ltwos{X\thetahat - X\thetastar}^2 =
\frac{2}{n} + \frac{1}{n} \sum_{i=2}^m (\sqrt{n}\thetahat_i -
\sqrt{n}(\thetastar_i - \thetastar_{m+i}))^2.
\end{align}
The first term on the right-hand side is obtained from the fact
$\ltwos{(X\thetahat - X\thetastar)_{1:2}}^2 = \ltwos{w_{1:2}}^2 = 2$,
recalling that the set-up assumes that the noise elements takes values
in $\{-1, 1\}$.

The right-hand side of equation~\eqref{eqn:simplified-loss-2} is a
Lasso objective function with design matrix $\sqrt{n}I_{(m-1)\times
  (m-1)}$. The second term on the right-hand side of
equation~\eqref{eqn:simplified-pred-error} is the associated
prediction error.  By choosing a proper $\lambda$ and using the fact
that $\thetastar$ is 2-sparse, it is well-known that the prediction
error scales as $\order(\frac{\log n}{n})$, which corresponds the fast
rate.  (Here we have recalled that the dimension of the Lasso problem
is $m - 1 = n - 2 $.)
%%%%%%%%%%%%%%%%%%%%%%%%%%%%%%%%%%%%%%%%%%%%%%%%%%%%%%%%%%%%%%%%%%%%%%%%%%

\section{Proof of Lemma~\ref{LEMMA:DECOMPOSED-TERM-LOWER-BOUND}}
\label{AppDecomposedOne}

Given our definition of $\Xmat$ in terms of the matrix $A \in \real^{2
  \times 2}$, it suffices to prove the two lower bounds
\begin{subequations}
\begin{align}
 \label{EqnFirstLemmaBoundOne} 
& \norms{A(\thetahat_\lambda)_{1:2} - A\opt_{1:2}}_2^2\\ &\qquad \geq
 \indicator\Big[\lambda\gamma_1 > 4 B (\sin^2(\alpha) \radius +
   \frac{\ltwos{w_{1:2}}}{\sqrt{\numobs}} )
   \Big]\sin^2(\alpha)(\radius -2B)_+^2 \quad \mbox{and},\nonumber \\
\label{EqnFirstLemmaBoundTwo}
& \norms{A(\thetahat_\lambda)_{2i-1:2i} - A\opt_{2i-1:2i}}_2^2\\
& \qquad \geq
\indicator\Big[0 \leq \wprime_i \leq B\Big] (\frac{B^2}{4} - \lambda
\gamma_1) \quad \mbox{for $i = 2, 3, \ldots,\numobs/2$.} \nonumber
\end{align}
\end{subequations}

In the proofs to follow, it is convenient to omit reference to the
index $i$.  In particular, viewing the index $i$ as fixed a priori, we
let $\solver$ and $\uvec^*$ be shorthand representations of the
sub-vectors $(\thetahat_\lambda)_{2i-1,2i}$, and
$\thetastar_{2i-1,2i}$, respectively. We introduce the normalized
noise $\altnoise \defeq w_{2i-1,2i}/\sqrt{n}$. By our construction of
the design matrix $X$ in terms of $A$, the vector $\thetahat_\lambda$
is a local minimizer of the objective function if and only if
$\solver$ is a local minimum of the following loss:
\begin{align*}
\ell(u;\lambda) \defeq \norms{Au - A\uvec^* - \altnoise}_2^2 + \lambda
\reg_{2i-1}(u_1) + \lambda \reg_{2i}(u_2),
\end{align*}
where this statement should hold for each $i \in [n/2]$.  Hence, it
suffices to find a local minimum of $\ell(u;\lambda)$ such that the
bounds~\eqref{EqnFirstLemmaBoundOne} and~\eqref{EqnFirstLemmaBoundTwo}
hold.

%%%%%%%%%%%%%%%%%%%%%%%%%%%%%%%%%%%%%%%%%%%%%%%%%%%%%%%%%%%%%%%%%%%%%%%

\subsection{Proof of inequality~\eqref{EqnFirstLemmaBoundOne}}

If $\lambda\gamma_1 \leq 4B(\sin^2(\alpha) \radius +
\ltwos{\altnoise})$, then the lower
bound~\eqref{EqnFirstLemmaBoundOne} is trivial, and in particular, it
will hold for
\[
	\solver \defeq \arg\min_{u} \ell(u;\lambda).
\]

Otherwise, we may assume that $\lambda\gamma_1 > 4B(\sin^2(\alpha)
\radius + \ltwos{\altnoise})$.  In this case ($i = 1$), we have $u^* =
(\radius/2, \radius/2)$.  Defining the vectors $v^* \defeq Au^* = (0,
\sin(\alpha) \radius)$ and $\util \defeq (0, 0)$, we have
\begin{align}
\label{eqn:statement-1-center-bound}
    \ell(\widetilde u;\lambda) &= \norms{A\widetilde u - v^* -
      \altnoise}_2^2 + \lambda \reg_1(\widetilde u_1) +
    \lambda\reg_2(\widetilde u_2) = \norms{v^* + \altnoise}_2^2.
\end{align}
We claim that
\begin{align}
\label{EqnUclaim}
\inf_{u\in \partial U} \ell(u;\lambda) > \ell(\widetilde u;\lambda),
\end{align}
where $U \defeq \{u \in \real^2 \, \mid \, |u_1| \leq B \mbox{ or }
|u_2| \leq B\}$, and $\partial U$ denotes its boundary.  If $\rho$ is
a convex function, then the lower bound~\eqref{EqnUclaim} implies that
any minimizers of the function $\ell(\cdot; \lambda)$ lie in the
interior of $U$. Otherwise, it implies that at least one local
minimum---say $\solver$---lies in the interior of $U$.  Since
$\solver_1 \leq B$ and $\solver_2 \leq B$, we have the lower bound
\begin{align*}
\norms{A(\thetahat_\lambda)_{1:2} - A\opt_{1:2}}_2^2 &=
\norms{A\solver - v^*}_2 = \cos^2(\alpha)(\solver_1 - \solver_2)^2 +
\sin^2(\alpha)(\radius - \solver_1 - \solver_2)^2\\ &\geq
\sin^2(\alpha)(\radius - \solver_1 - \solver_2)^2 \geq
\sin^2(\alpha)(\radius - 2B)_+^2,
\end{align*}
which completes the proof. \\

It remains to prove the lower bound~\eqref{EqnUclaim}.  For any $u\in
\partial U$, we have
\begin{align}
\label{eqn:statement-1-boundary-bound}
\ell(u;\lambda) &= \norms{Au - v^* - \altnoise}_2^2 + \lambda
\reg_1(\widetilde u_1) + \lambda\reg_2(\widetilde u_2)
\stackrel{(i)}{\geq} \norms{Au + v^* + \altnoise}_2^2 +
\lambda\gamma_1 \nonumber\\ &\stackrel{(ii)}{\geq} \norms{v^* +
  \altnoise}_2^2 + 2(v^* + \altnoise)^T A u + \lambda\gamma_1.
\end{align}
Inequality (i) holds since either $\widetilde u_1$ or $\widetilde u_2$
is equal to $B$, and \mbox{$\min\{ \reg_1(B),\reg_2(B) \} \geq
  \gamma_1$} by definition, whereas inequality (ii) holds since
\mbox{$\ltwos{a+b}^2 \geq \ltwos{b}^2 + 2 b^Ta$.} We notice that
\begin{align*}
& \inf_{u\in \partial U} 2(v^* + \altnoise)^T A u \geq \inf_{u\in
    \partial U} \left\{ 2\langle v^*, A u \rangle - 2
  \ltwos{\altnoise} \ltwo{ A u} \right\}\\ & \qquad = \inf_{u\in
    \partial U} \left\{ 2 \sin^2(\alpha) \radius(u_1+u_2) - 2
  \ltwos{\altnoise}\sqrt{\cos^2(\alpha) (u_1+u_2)^2 + \sin^2(\alpha)
    (u_1-u_2)^2}\right\} \\ & \qquad \geq - 4B(\sin^2(\alpha) \radius
  + \ltwos{\altnoise}).
\end{align*}
Combining this lower bound with
inequality~\eqref{eqn:statement-1-boundary-bound} and the bound
$\lambda\gamma_1 > 4B(\sin^2(\alpha) \radius + \ltwos{\altnoise})$
yields the claim~\eqref{EqnUclaim}.

%%%%%%%%%%%%%%%%%%%%%%%%%%%%%%%%%%%%%%%%%%%%%%%%%%%%%%%%%%%%%%%%%%%%%

\subsection{Proof of inequality~\eqref{EqnFirstLemmaBoundTwo}}

For $i=2,3,\dots,\numobs/2$, consider $u^* = (0,0)$ and recall the
vector $a_i = (\cos \alpha, \sin \alpha)$, as well as our assumption
that
\begin{align*}
\gamma_i \defeq \reg_{2i-1}(B)/B \leq \frac{\reg_{2i}(B)}{B}.
\end{align*}
Define the vector $\widetilde u \defeq (a_i^T\altnoise, 0)$ and let
\mbox{$\solver = \arg \min_{u \in \real^2} \ell(u;\lambda)$} be an
arbitrary global minimizer.  We then have
\begin{align*}
\norms{A\solver - \altnoise}_2^2 
& \leq \norms{A\solver -
  \altnoise}_2^2 + \lambda \reg_1(\solver_1) + \lambda
\reg_2(\solver_2) \leq \norms{A\widetilde u - \altnoise}_2^2 + \lambda
\reg_1(\widetilde u_1) + \lambda \reg_2(\widetilde u_2),
\end{align*}
since the regularizer is non-negative, and $\solver$ is a global minimum.
Using the definition of $\util$, we find that
\begin{align*}
\norms{A\solver - \altnoise}_2^2 & \leq \ltwos{a_i a_i^T\altnoise -
  \altnoise}^2 + \lambda\reg_1(a_i^T\altnoise) = \norms{\altnoise}_2^2 -
(a_i^T\altnoise)^2 + \lambda\reg_1(a_i^T\altnoise),
\end{align*}
where the final equality holds since $a_i a_i^T$ defines an orthogonal
projection. By the triangle inequality, we find that $\norms{A\solver
  - \altnoise}_2^2 \geq \norms{\altnoise}_2^2 - \norms{A\solver}_2^2$,
and combining with the previous inequality yields
\begin{align}
\label{eqn:error-with-zero-coefficient-1}
\norms{A\solver}_2^2 \geq (a_i^T\altnoise)^2 - \lambda
\reg_1(a_i^T\altnoise).
\end{align}
Now if $B/2 \leq a^T\altnoise \leq B$, then we have
$\reg_1(a_i^T\altnoise) \leq \reg_1(B) = \gamma_i \leq \gamma_1$.
Substituting this relation into
inequality~\eqref{eqn:error-with-zero-coefficient-1}, we have
$\norms{A\solver - \altnoise}_2^2 \geq \indicator(B/2 \leq
a_i^T\altnoise\leq B) \left( B^2/4 - \lambda\gamma_1\right)$, which
completes the proof.

%%%%%%%%%%%%%%%%%%%%%%%%%%%%%%%%%%%%%%%%%%%%%%%%%%%%%%%%%%%

\section{Proof of Lemma~\ref{LEMMA:GRADIENT-DESCENT-TERMWISE-LOWER-BOUND}}
\label{AppDecomposedTwo}

Similar to the proof of Lemma~\ref{LEMMA:DECOMPOSED-TERM-LOWER-BOUND},
it is convenient to omit reference to the index $i$.  We let $\uvec^t$
and $\uvec^*$ be shorthand representations of the sub-vectors
$\uvec^t_{2i-1,2i}$, and $\thetastar_{2i-1,2i}$, respectively. We
introduce the normalized noise $\altnoise \defeq
w_{2i-1,2i}/\sqrt{n}$. By our construction of the design matrix $X$
and the update formula~\eqref{eqn:local-descent-update-formula}, the
vector $\uvec^t$ satisfies the recursion
\begin{subequations}
\begin{align}
\label{eqn:u-gradient-update-formula}
  \uvec^{t+1} \in \underset{\ltwos{u - \uvec^t}\leq \beta}{\argmin}
  \ell(u;\lambda),
\end{align}
where $\beta \defeq \ltwos{u^{t+1} - \uvec^t} \leq \eta$ and
the loss function takes the form
\begin{align}
\label{eqn:define-square-norm-T}
\ell(u;\lambda) \defeq \underbrace{\norms{A\uvec - A\uvec^* -
    \altnoise}_2^2}_{\defeq T} + \lambda \reg_{2i-1}(\uvec_1) +
\lambda \reg_{2i}(\uvec_2).
\end{align}
\end{subequations}
This statement holds for each $i \in [n/2]$.  Hence, it suffices to
study the update formula~\eqref{eqn:u-gradient-update-formula}.

%%%%%%%%%%%%%%%%%%%%%%%%%%%%%%%%%%%%%%%%%%%%%%%%%%%%%%%%%%%%%%%%%%%%%

\subsection{Proof of part (a)}  

For the case of $i=1$, we assume that the event $\event_2$
holds. Consequently, we have $\max\{\uvec^0_1,\uvec^0_2\}\leq 0$ and
$\lambda\gamma_1 \geq 2\sin^2(\alpha)\altradius
+2\ltwos{\altnoise}$. The corresponding regression vector is $\ustar =
(\altradius/2, \altradius/2)$. Let us define
\begin{align*}
b_1 \in \arg\sup_{u\in(0,B]} \rho_{2i-1}(u) \quad \mbox{and} \quad b_2
  \in \arg\sup_{u\in(0,B]} \rho_{2i}(u).
\end{align*}
Our assumption implies \mbox{$\rho'_{2i-1}(b_1) \geq \gamma_1$} and
\mbox{$\rho'_{2i}(b_2) \geq \gamma_1$.}  We claim that
\begin{align}
\label{claim:i1-bound}
 \uvec^t_k \leq b_k \leq B \quad \mbox{for $k=1,2$ and for all
   iterations $t = 0, 1, 2, \ldots$.}
\end{align}
If the claim is true, then we have
\begin{align*}
\norms{A\theta^t_{1:2} - A\opt_{1:2}}_2^2 & = \cos^2(\alpha)(\uvec^t_1
- \uvec^t_2)^2 + \sin^2(\alpha)(\altradius - \uvec^t_1 - \uvec^t_2)^2
\\
& \geq \sin^2(\alpha)(\altradius - \uvec^t_1 - \uvec^t_2)^2 \geq
\sin^2(\alpha)(\altradius - 2B)_+^2 \\
& \geq \frac{\sigma \altradius}{4\sqrt{n}},
\end{align*}
where the final inequality follows from substituting the definition of
$\alpha$, and using the fact that $2B \leq r/2$. Thus, it suffices to
prove the claim~\eqref{claim:i1-bound}.

We prove the claim~\eqref{claim:i1-bound} by induction on the
iteration number $t$.  It is clearly true for $t=0$. Assume that the
claim is true for a specific integer $t \geq 0$, we establish it for
integer $t+1$. Our strategy is as follows: suppose that the vector
$\uvec^{t+1}$ minimizes the function $\ell(\uvec;\lambda)$ inside the
ball $\{\uvec:\ltwos{\uvec - \uvec^t}\leq \beta\}$.  Then, the scalar
$\uvec^{t+1}_1$ satisfies
\begin{align*}
\uvec^{t+1}_1 = \underset{x :\ltwos{(x,u^{t+1}_2)-\uvec^t}\leq
  \beta}{\argmin} f(x) \quad \mbox{where}~ f(x) \defeq
\ell((x,\uvec^{t+1}_2);\lambda).
\end{align*}
We now calculate the generalized derivative~\cite{Clarke83} of the
function $f$ at $\uvec^{t+1}_1$. It turns out that
\begin{align}
\label{eqn:claim-i1-bound-basic-stat}
\mbox{either}\quad \uvec^{t+1}_1 \leq \uvec^t_1 \leq b_1,\quad
\mbox{or}\quad \partial f(\uvec^{t+1}_1)\cap(-\infty,0] \neq
  \emptyset.
\end{align}
Otherwise, there is a sufficiently small scalar $\delta > 0$ such that
\begin{align*}
\ltwos{(u^{t+1}_1-\delta,u^{t+1}_2)-\uvec^t}\leq \beta
\quad\mbox{and}\quad f(u^{t+1}_1-\delta) < f(u^{t+1}_1),
\end{align*}
contradicting the fact that $\uvec^{t+1}_1$ is the minimum point.  In
statement~\eqref{eqn:claim-i1-bound-basic-stat}, if the first
condition is true, then we have $\uvec^{t+1}_1 \leq b_1$. We claim
that the second condition also implies $\uvec^{t+1}_1 \leq b_1$.

To prove the claim, we assume by contradiction that $\uvec^{t+1}_1 > b_1$ and $\partial f(\uvec^{t+1}_1)\cap(-\infty,0] \neq \emptyset$.
Note that the function $f(x)$ is differentiable for $x > 0$. In particular, for $\uvec^{t+1}_1 > b_1$, we have
\begin{align}\label{eqn:equation-partial-T-u1}
  f'(\uvec^{t+1}_1) &= \frac{\partial T}{\partial \uvec_1} \Big|_{u_1 = \uvec^{t+1}_1} + \lambda \rho_{2i-1}'(\uvec^{t+1}_1)\nonumber\\
  &= - 2( \sin^2(\alpha)\altradius + a_i^T \altnoise) + 2\uvec^{t+1}_1 - 2(1 - 2\sin^2(\alpha))\uvec^{t+1}_2 + \lambda \rho_{2i-1}'(\uvec^{t+1}_1).
\end{align}
Since $1 - 2\sin^2(\alpha) \leq 1$ and $a_i^T \altnoise \leq \ltwos{\altnoise}$, and also because $\uvec^{t+1}_1 > b_1$ and $\uvec^{t+1}_2 \leq \uvec^{t}_2 + \beta \leq b_2 + \beta$, equation~\eqref{eqn:equation-partial-T-u1} implies
\begin{align}
\label{eqn:neg-gradient-least-square-loss}
f'(\uvec^{t+1}_1) \geq -2\sin^2(\alpha)\altradius - 2\ltwos{\altnoise} + 2(b_1 - b_2 - \beta) + \lambda \rho_{2i-1}'(\uvec^{t+1}_1).
\end{align}
Recall that $b_1,b_2\in [0,B]$, and also using the fact that
\begin{align*}
\reg'_{1}(\uvec^{t+1}_1) \geq \reg'_{1}(b_1) - \beta\curbound \geq
\gamma_1 - \beta\curbound,
\end{align*}
we find that
\begin{align}
f'(\uvec^{t+1}_1) &\geq -2\sin^2(\alpha)\altradius - 2
\ltwos{\altnoise} - 2(B+\beta) + \lambda(\gamma_1 -
\beta\curbound).\nonumber \\ 
 \label{eqn:neg-gradient-loss}
& \geq -2\sin^2(\alpha)\altradius - 2 \ltwos{\altnoise} - 3B +
 \lambda\gamma_1.
\end{align}
Here the second inequality follows since $\beta \leq \eta \leq
\min\{B, \frac{B}{\lambda\curbound}\}$.  Since the inequality
$\lambda\gamma_1 > 2\sin^2(\alpha)\altradius + 2\ltwos{\altnoise} +
3B$ holds, inequality~\eqref{eqn:neg-gradient-loss} implies that
\mbox{$f'(\uvec^{t+1}_1) > 0$.}  But this conclusions contradicts the
assumption that
\begin{align*}
\partial f(\uvec^{t+1}_1) \cap (-\infty,0] \neq \emptyset.
\end{align*}
Thus, in both cases, we have $\uvec^{t+1}_1 \leq b_1$.

The upper bound for $\uvec^{t+1}_2$ can be proved following the same
argument.  Thus, we have completed the induction.

%%%%%%%%%%%%%%%%%%%%%%%%%%%%%%%%%%%%%%%%%%%%%%%%%%%%%%%%%%%%%%%%%%%%%%%%%%%

\subsection{Proof of part (b)}  

For the case of $i=2, 3, \ldots, \numobs/2$, we assume that the event
$i\in \Sset_2$ holds.  Consequently, we have $\lambda\gamma_1 < 2
a_i^T \altnoise - 4 B$ as well as our assumption that
\begin{align*}
\sup_{u\in(0,B]} \rho_{2i-1}'(u) = \gamma_i \leq \gamma_1.
\end{align*}
The corresponding regression vector is $\ustar = (0,0)$.  Let $\uhack$
be the stationary point to which the sequence
$\{\uvec^t\}^\infty_{t=0}$ converges.  We claim that
\begin{align}
\label{eqn:claim-part-b}
|\uhack_1| \geq B \mbox{ or } |\uhack_2| \geq B.
\end{align}
If the claim is true, then by inequality $\cos^2(\alpha) \geq
\sin^2(\alpha)$ obtained from the definition of $\alpha$, we have
\begin{align*}
\norms{A\thetahack_{2i-1:2i} - A\opt_{2i-1:2i}}_2^2 &=
\cos^2(\alpha)(\thetahack_{2i-1} - \thetahack_{2i})^2 +
\sin^2(\alpha)(\thetahack_{2i-1} + \thetahack_{2i})^2\\ &\geq
2\sin^2(\alpha)((\uhack_1)^2 + (\uhack_2)^2) \geq 2\sin^2(\alpha) B^2
\\
& = \frac{\sigma^3/r}{8\numobs^{3/2}} \geq \frac{\sigma
  \altradius}{8\numobs^{3/2}}.
\end{align*}
which completes the proof of part (b). Thus, it suffices to proof the
claim~\eqref{eqn:claim-part-b}. \\

To prove the claim~\eqref{eqn:claim-part-b}, we notice that $\uhack$
is a local minimum of the loss function $\ell(\cdot;\lambda)$. Let
\begin{align*}
f(x) \defeq \ell((x,\uhack_2);\lambda)
\end{align*}
be a function that restricts the second argument of the function
$\ell(\cdot;\lambda)$ to be $\uhack_2$.  Since $\uhack$ is a local
minimum of the function $\ell(\cdot;\lambda)$, the scalar $\uhack_1$
must be a local minimum of $f$. Consequently, the zero vector must
belong to the generalized derivative~\cite{Clarke83} of $f$, which we
write as $0\in \partial f(\uhack_1)$. We use this fact to prove the
claim~\eqref{eqn:claim-part-b}.

Assume by contradiction that the claim~\eqref{eqn:claim-part-b} does
not hold, which means that $|\uhack_1| < B$ and $|\uhack_2| <
B$. Calculating the generalized derivative of $f$, we have
\begin{align}
\label{eqn:gradient-lemma-part-2-derivatives-1}
\partial f(\uhack_1) &= \{- 2 a_i^T \altnoise + 2 \uhack_1 - 2(1 -
2\sin^2(\alpha))\uhack_2 + \lambda g_1:~ g_1 \in \partial
\reg_{2i-1}(\uhack_1) \}.
\end{align}
Using the fact that $|\uhack_1| < B$, $|\uhack_2| < B$ and the upper
bound
\begin{align*}
|g_1| \leq \sup_{u\in(0,B]}\reg'_{2i-1}(u) = \gamma_i \leq \gamma_1,
\end{align*}
the derivative~\eqref{eqn:gradient-lemma-part-2-derivatives-1} is
upper bounded by
\begin{align*}
- 2 a_i^T \altnoise + 2 \uhack_1 - 2(1 - 2\sin^2(\alpha))\uhack_2 +
\lambda g_1 < - 2 a_i^T \altnoise + 4 B + \lambda \gamma_1.
\end{align*}
Under the assumption that $\lambda \gamma_1 < 2 a_i^T \altnoise - 4
B$, the above inequality implies $0\notin \partial f(\uhack_1)$,
contradicting the fact that $0\in \partial f(\uhack_1)$.

%%%%%%%%%%%%%%%%%%%%%%%%%%%%%%%%%%%%%%%%%%%%%%%%%%%

\subsection{Proof of part (c)}

Recall that the $2$-vector $\theta^0_{1:2}$ follows a
$N(0,b^2I_{2\times 2})$ distribution, so that
\begin{align*}
\mprob[\event_0] = \mprob\big[\max\{\theta^0_1,\theta^0_2\}\leq 0\big]
= \frac{1}{4},
\end{align*}
which establishes the first claim of part (c). To prove the second
statement, we notice that the inequality $2\sin^2(\alpha)\altradius +
\frac{2 \ltwos{w_{1:2}}}{\sqrt{\numobs}} + 3B \leq 2\wprime_i - 4 B$
can be written in the equivalent form
\begin{align}
\label{eqn:second-event-s3}
2a_i^T w_{2i-1:2i} - 2\ltwos{w_{1:2}} \geq \frac{7\sigma}{4}.
\end{align}
The inequality~\eqref{eqn:second-event-s3} holds if
$a_i^Tw_{2i-1:2i}/\sigma \geq 1$ and $\ltwos{w_{1:2}}^2/\sigma^2 \leq
1/64$. In fact, $a_i^T w_{2i-1:2i}/\sigma$ is distributed as $N(0,1)$;
and $\ltwos{w_{1:2}}^2/\sigma^2$ follows a chi-square distribution
with two degrees of freedom. The two events are independent, and each
of them happens with (positive) constant probability. Thus, there is a
numerical constant $c > 0$ such that
inequality~\eqref{eqn:second-event-s3} holds with probability at
least~$c$.

%%%%%%%%%%%%%%%%%%%%%%%%%%%%%%%%%%%%%%%%%%%%%%%%%%%%%%%%%%%%%%%%%%%%%%%%

\bibliographystyle{abbrvnat} \bibliography{bib}

\end{document}